\documentclass[preprint,12pt]{elsarticle}

\usepackage{graphicx}
\usepackage{color}
\usepackage{caption}
\usepackage{subcaption}
\usepackage{amsmath,amssymb,amsfonts,bm}
\usepackage{lineno}

\usepackage[normalem]{ulem}
\usepackage{linenofix}
\DeclareMathOperator*{\argmin}{arg\,min}

\begin{document}

\begin{frontmatter}

%% Title, authors and addresses

\title{Learning Parameters and Constitutive Relationships with Physics Informed Deep Neural Networks}

%% use the tnoteref command within \title for footnotes;
%% use the tnotetext command for the associated footnote;
%% use the fnref command within \author or \address for footnotes;
%% use the fntext command for the associated footnote;
%% use the corref command within \author for corresponding author footnotes;
%% use the cortext command for the associated footnote;
%% use the ead command for the email address,
%% and the form \ead[url] for the home page:
%%
%% \title{Title\tnoteref{label1}}
%% \tnotetext[label1]{}
%% \author{Name\corref{cor1}\fnref{label2}}
%% \ead{email address}
%% \ead[url]{home page}
%% \fntext[label2]{}
%% \cortext[cor1]{}
%% \address{Address\fnref{label3}}
%% \fntext[label3]{}

%% use optional labels to link authors explicitly to addresses:
%% \author[label1,label2]{<author name>}
%% \address[label1]{<address>}
%% \address[label2]{<address>}

\author[pnnl]{Alexandre M. Tartakovsky\corref{cor1}}
\author[pnnl]{Carlos Ortiz Marrero}
\author[upenn]{Paris Perdikaris}
\author[pnnl]{Guzel D. Tartakovsky}
\author[pnnl]{David Barajas-Solano}

\address[pnnl]{Pacific Northwest National Laboratory}
\address[upenn]{University of Pennsylvania}
\cortext[cor1]{alexandre.tartakovsky@pnnl.gov}

\begin{abstract}
We present a physics informed deep neural network (DNN) method for estimating parameters and unknown physics (constitutive relationships) in partial differential equation (PDE) models. We use PDEs in addition to measurements to train DNNs to approximate unknown parameters and constitutive relationships as well as states. The proposed  approach increases the accuracy of DNN approximations of partially known functions when a limited number of measurements is available and allows for training DNNs when no direct measurements of the functions of interest are available. We employ physics informed DNNs to estimate the unknown space-dependent diffusion coefficient in a linear diffusion equation and an unknown constitutive relationship in a non-linear diffusion equation. For the parameter estimation problem, we assume that partial measurements of the coefficient and states are available and demonstrate that under these conditions, the proposed method is more accurate than state-of-the-art methods. For the non-linear diffusion PDE model with a fully unknown constitutive relationship (i.e., no measurements of constitutive relationship are available), the physics informed DNN method can accurately estimate the non-linear constitutive relationship based on state measurements only.  Finally, we demonstrate that the proposed method remains accurate in the presence of measurement noise. 
\end{abstract}

\begin{keyword}
Deep Neural Networks \sep partial differential equations \sep parameter estimation \sep learning unknown physics
%% keywords here, in the form: keyword \sep keyword

%% MSC codes here, in the form: \MSC code \sep code
%% or \MSC[2008] code \sep code (2000 is the default)

\end{keyword}

\end{frontmatter}

%%
%% Start line numbering here if you want
%%
%\linenumbers

\section{Introduction}

Physical models of many complex natural systems are, at best, ``partially'' known as conservation laws do not provide a closed system of equations. Accurate theoretical models for closing the system of conservation equations are available for homogeneous systems exhibiting time and length scale separation. Examples of accurate closures include Newtonian stress for homogeneous (Newtonian) fluids, Fick's law for mass flux in diffusion processes, and the Darcy law for fluid flux in porous media. For more complex systems, including non-homogeneous turbulence, non-Newtonian fluid flow, multiphase flow and transport in porous media, and granular materials, accurate theoretical closures are not available. Instead,  phenomenological constitutive relationships are used, which are usually accurate for a narrow range of conditions. Even when sufficiently accurate closed-form partial differential equations (PDE) models are available, (space-dependent) parameters are typically unknown.

Computational approaches for parameter and constitutive law estimation cast this inherently ill-posed problem as an optimization or a statistical inference task, usually requiring  repeated evaluation of expensive forward solvers until the parameters that minimize a given error metric are found or until space of parameter configurations satisfying available data is explored. This 
  results in significant and often intractable computational cost.
  Furthermore, minimization via gradient-based methods requires computing gradients from said expensive forward models, which either requires additional computational cost or careful formulation of the adjoint problems.
  Additionally, forward modeling requires the knowledge of initial and boundary conditions, which usually are not fully known and, therefore, also must be estimated from data together with unknown parameters and constitutive relations, significantly complicating parameter estimation.

Although significant progress has been made over the last two decades involving high-order schemes for PDEs, automatic differentiation of computer code, PDE-constrained optimization, and optimization under uncertainty, parameter estimation in large-scale problems remains a significant challenge \cite{lieberman2010parameter}.

While there are a number of established methods for parameter estimation in (closed-form) PDE models,
  such as Bayesian inference and maximum a posteriori probability (MAP) estimation~\cite{stuart2010inverse,carrera2005inverse}, %
existing approaches for learning unknown physics from partially known models and data are few and not fully mature.

Unknown forms of equations at a given scale can either be found by upscaling (coarsening) known equations governing the same process at a smaller scale (e.g., \cite{whitaker2013method}), or learned from data. In this work, we are interested in the latter.

In the most general case, system dynamics can be described as 
\begin{equation*}
  \frac{\partial u}{\partial t} = F(u),
\end{equation*}
where $F$ is a function or differential operator. 
A recent review of methods for learning $F$ can be found in \cite{rudy2018data}. These methods include NARMAX \cite{chen1989representations}, equation-free methods \cite{kevrekidis2003equation}, Laplacian spectral analysis \cite{giannakis2012nonlinear}, and neural networks \cite{gonzalez1998identification}. 
Recently, data-driven approximations of Koopman operators, including dynamic mode decomposition \cite{williams2015data}, diffusion maps \cite{williams2015data}, delay coordinates \cite{giannakis2017data,brunton2017chaos,askham2018variable}, and neural networks \cite{yeung2017learning,wehmeyer2018time, mardt2018vampnets,lusch2017deep,raissi2018deep}, have gained significant attention.
This work concerns with problems where conservation laws and other physical knowledge provide significantly more information
about the structure of the operator $F(u)$.
Specifically, we are interested in a steady state of the problem
\begin{equation}
  \label{eq:pde-problem}
\frac{\partial u(\mathbf{x},t)}{\partial t} = -\nabla \cdot (K(\mathbf{x},u) \nabla u(\mathbf{x},t)),
\end{equation}
where $K(\mathbf{x},u)$
is an unknown constitutive relationship, which is a function of space and the PDE state. The boundary conditions may or may not be known.  
Among other physical phenomena, this equation describes flow in porous media~\cite{bear2013dynamics}.
  Two special cases of this problem are $K(\mathbf{x}, u) = K(\mathbf{x})$, where Eq (\ref{eq:pde-problem}) becomes a linear diffusion equation with heterogeneous diffusion coefficient $K(\mathbf{x})$; and $K(\mathbf{x}, u) = K(u)$, where Eq (\ref{eq:pde-problem}) becomes a nonlinear diffusion equation with a state-dependent coefficient $K(u)$.
We are interested in learning $K(\mathbf{x},u)=K(\mathbf{x})$ when measurements of both $u$ and $K$ are available, and $K(\mathbf{x},u)=K(u)$ when only $u$ measurements are available.

To achieve this objective, we propose a physics informed DNN method, where PDEs
and data are used to train DNN representations of the PDE states and the unknown parameters constitutive relations. 
With application to Eq (\ref{eq:pde-problem}), this approach consists of defining two DNNs, one for $K(\mathbf{x},u)$ and another for $u(\mathbf{x})$, together with auxiliary DNNs obtained by substituting $K$ and $u$ into Eq (\ref{eq:pde-problem}) and using automatic differentiation to evaluate the right-hand side expression and boundary conditions.
These networks are trained simultaneously using available data.
Our work extends the physics informed DNN method proposed in~\cite{raissi2017physicsPart1,raissi2017physicsPart2} for finding unknown constants and solutions to PDEs given states observations. It is important to note that the physics informed DNN method was originally developed for time-dependent PDEs, and the extension of this method proposed herein also can be applied to time-dependent problems. 

The work is organized as follows: in Section \ref{philms}, we introduce the PDE model and the physics-informed DNN approach.
In Section \ref{linearPDE}, we provide a detailed study of the accuracy and performance of the proposed approach for a linear parameter estimation problem with $K(\mathbf{x},u)=K(\mathbf{x})$.
In Section \ref{nonlinearPDE}, we demonstrate the method's accuracy for learning the non-linear constitutive relationship $K(\mathbf{x},u) = K(u)$.
Discussion and conclusions are presented in Section \ref{conclusions}.    

\section{Physics-informed Deep Neural Network Approach}\label{philms}

Consider the steady-state PDE
\begin{equation}
  \label{eq:steady-pde}
  \mathcal{L}[u,K(\mathbf{x},u)] = 0, \quad \mathbf{x} \in \Omega
\end{equation}
subject to the Dirichlet and Neumann boundary conditions
\begin{align}
  \label{eq:steady-pde-bc-D}
  u(\mathbf{x}) &=g(\mathbf{x}), && \mathbf{x}\in \partial_D \Omega,\\
  \label{eq:steady-pde-bc-N}
  \mathbf{n} \cdot K(\mathbf{x},u) \nabla u(\mathbf{x}) &= q(\mathbf{x}), && \mathbf{x}\in \partial_N \Omega.
\end{align}
Here, $\mathcal{L}$ is the known differential operator, %
  $\Omega \subset \mathbb{R}^d$, $d \in [1, 3]$ is the simulation domain with the boundary $\partial \Omega$, $\partial_D \Omega$ and $\partial_N \Omega$ are the ``Dirichlet'' and ``Neumann'' parts of the boundary satisfying $\partial \Omega = \partial_D \Omega \cup \partial_N \Omega$ and $\partial_D \Omega \cap \partial_N \Omega \equiv \varnothing$, and $\mathbf{n}$ is the outward unit vector normal to $\partial \Omega$.
  We denote the state of the model by $u : \Omega \to \mathcal{U} \subseteq \mathbb{R}$ and the unknown constitutive relation by $K : \Omega \times \mathcal{U} \to \mathcal{K} \subseteq \mathbb{R}$.

We assume that $N_K$ measurements of $K$, $N_u$ measurements of $u$, $N_D$ measurements of $g$, and $N_N$ measurements of $q$ are 
  collected at the locations $\{ \mathbf{x}^K_i \}^{N_K}_{i = 1}$, $\{ \mathbf{x}^u_i \}^{N_u}_{i = 1}$, $\{ \mathbf{x}^D_i \}^{N_D}_{i = 1}$, and $\{ \mathbf{x}^N_i \}^{N_N}_{i = 1}$, respectively.
  The observations are denoted by $K^{*}_i \equiv K(\mathbf{x}^K_i, u(\mathbf{x}^K_i))$ ($i = 1, \dots, N_K$), $u^{*} \equiv u(\mathbf{x}^u_i)$ ($i = 1, \dots, N_u$), $g^{*}_i \equiv g(\mathbf{x}^D_i)$ ($i = 1, \dots, N_D$), and $q^{*}_i \equiv q(\mathbf{x}^N_i)$ ($i = 1, \dots, N_N$).

To learn $K(\mathbf{x},u)$, we define the following DNNs for $u(\mathbf{x})$ and $K(\mathbf{x},u)$:
\begin{equation}
  \hat{u}(\mathbf{x};\theta) = \mathcal{NN}_u(\mathbf{x};\theta), \quad \hat{K}(\mathbf{x}, u; \gamma) = \mathcal{NN}_K(\mathbf{x}, u; \gamma),
  \label{DNNs}
\end{equation}
where $\theta$ and $\gamma$ are the DNN parameters.
  Substituting these two DNNs into the governing equation (\ref{eq:steady-pde}) and the Neumann boundary condition (\ref{eq:steady-pde-bc-N}), and evaluating the corresponding spatial derivatives via automatic differentiation, %
yields two additional ``auxiliary'' DNNs:
\begin{gather*}
  f(\mathbf{x}; \gamma, \theta) = \mathcal{L} \left [ \mathcal{NN}_u(\mathbf{x}; \theta),  \mathcal{N}\mathcal{N}_K \left ( \mathbf{x}, \mathcal{NN}_u \left (\mathbf{x}; \theta \right ); \gamma \right )  \right ] = \mathcal{NN}_f (\mathbf{x};\theta,\gamma),\\
  f_N(\mathbf{x}; \gamma, \theta) = \mathbf{n} \cdot \mathcal{NN}_K \left ( \mathbf{x},   \mathcal{NN}_u(\mathbf{x}; \theta); \gamma \right ) \nabla \mathcal{NN}_u(\mathbf{x};\theta) = \mathcal{N} \mathcal{N}_N (\mathbf{x};\theta,\gamma).
\end{gather*}
Next, we define the following loss function to train these four networks simultaneously:
\begin{equation}
  \label{loss_fn_with_BC} 
  \begin{split}
    L(\theta, \gamma) &= \frac{1}{N_K}\sum \limits_{i=1}^{N_K} \left [ \hat{K}(\mathbf{x}^K_i, \hat{u}(\mathbf{x}^K_i;\theta); \gamma) - K^*_i \right]^2 + \frac{1}{N_u}\sum \limits_{i=1}^{N_u} \left [\hat{u}(\mathbf{x}^u_i;\theta) - u^*_i \right]^2 \\
    &+ \frac{1}{N_D}\sum \limits_{i=1}^{N_D} \left [ \hat{u}(\mathbf{x}^D_i; \theta) - g^{*}_i \right ]^2 +\frac{1}{N_N}\sum \limits_{i=1}^{N_N} \left [ f_N(\mathbf{x}^N_i;\gamma,\theta) - q^{*}_i \right ]^2\\
    &+ \frac{1}{N_c} \sum \limits_{i=1}^{N_c} f(\mathbf{x}^c_i;\gamma,\theta)^2.
  \end{split}
\end{equation}
The first and second terms in $L$ force the $K$ and $u$ DNNs to match the $K$ and $u$ measurements.
  The third and fourth terms enforce Dirichlet and Neumann boundary conditions.
  Finally, the fifth term enforces the PDE at $N_c$ ``collocation'' points $\{ \mathbf{x}^c_i \}^{N_c}_{i = 1}$ that can be chosen uniformly or non-uniformly over $\Omega$ depending on the problem.

The DNNs are trained, i.e., $\theta$ and $\gamma$ are found, by minimizing the loss function:
\begin{equation}
  \label{eq:min_loss}
  (\theta, \gamma) = \argmin_{\theta,\gamma} L(\theta, \gamma).
\end{equation}
The minimization is carried out using the L-BFGS-B method~\cite{byrd1995limited} together with Xavier's normal initialization scheme~\cite{glorot2010understanding}.
  We use a quasi-Newton optimizer such as L-BFGS-B instead of stochastic gradient descent \cite{zhang2013asynchronous} (a more common optimizer for DNNs) because of its superior rate of convergence and more favorable computational cost for problems with a relatively small number of observed data.

Note that the proposed method does not make any assumptions about the measurement noise.  Also, our method can be easily extended to time-dependent PDEs by defining DNNs in Eq (\ref{DNNs}) as functions of both $\mathbf{x}$ and $t$ \cite{raissi2017physicsPart1,raissi2017physicsPart2}.

In the following two sections, we apply the physics-informed DNNs to learn parameters and constitutive relationships in PDE models
of the form (\ref{eq:steady-pde})--(\ref{eq:steady-pde-bc-N}).

\section{Parameter estimation in a linear diffusion equation}
\label{linearPDE}

In this section, we consider a linear diffusion equation with unknown diffusion coefficient $K(\mathbf{x})$,
\begin{equation}
  \label{eq:linearPDE}
  \nabla \cdot (K(\mathbf{x}) \nabla u(\mathbf{x})) = 0,  \quad \mathbf{x} \equiv (x_1,x_2)^T \in (0,1) \times (0,1)
\end{equation}
subject to the Dirichlet boundary conditions
\begin{equation}
  \label{eq:linearPDE-bc-D}
  u(\mathbf{x}) = 1, \quad x_2 = 0 \quad \textnormal{and}\quad 
  u(\mathbf{x}) = 0, \quad x_2 = 1
\end{equation}
and the Neumann boundary conditions
\begin{equation}
  \label{eq:linearPDE-bc-N}
  \frac{\partial u(\mathbf{x})}{\partial x_1} = 0 \quad x_1 = \{ 0, 1 \}.
\end{equation}
Among other problems, this equation describes saturated flow in heterogeneous porous media with hydraulic conductivity $K(\mathbf{x})$~\cite{bear2013dynamics}.
We assume that $N_K$ measurements of $K(\mathbf{x})$ and $N_u$ measurements of $u(\mathbf{x})$ are available:
$K^{*}_i \equiv K(\mathbf{x}^K_i)$ $(i = 1, \dots, N_K)$ and $u^{*}_i \equiv u(\mathbf{x}^u_i)$ $(i = 1, \dots, N_u)$.
We define DNNs for $K(\mathbf{x})$ and $u(\mathbf{x})$, $ \hat{K}(\mathbf{x};\gamma)=\mathcal{N} \mathcal{N}_K(\mathbf{x};\gamma) $, and $ \hat{u}(\mathbf{x};\theta)=\mathcal{N} \mathcal{N}_u(\mathbf{x};\theta), $ together with two auxiliary DNNs $ f(\mathbf{x};\gamma,\theta)=\nabla \cdot [\mathcal{N} \mathcal{N}_K(\mathbf{x};\gamma) \nabla \mathcal{N} \mathcal{N}_u(\mathbf{x};\theta)] = \mathcal{N} \mathcal{N}_f (\mathbf{x};\theta,\gamma) $ and $f_N(\mathbf{x};\theta)= \partial \mathcal{N} \mathcal{N}_u(\mathbf{x};\theta)/\partial x_2 = \mathcal{N} \mathcal{N}_N (\mathbf{x};\theta)$.
For this problem, the loss function takes the form:
\begin{equation}
  \label{eq:loss_fn_linear} 
  \begin{split}
    L(\theta, \gamma) &= \frac{1}{N_K} \sum \limits_{i=1}^{N_K} \left [ \hat{K}(\mathbf{x}^K_i;\gamma) - K^{*}_i \right ]^2 + \frac{1}{N_u} \sum \limits_{i=1}^{N_u} \left [ \hat{u}(\mathbf{x}^u_i;\theta) - u^{*}_i \right ]^2 \\ 
        &+ \frac{1}{N_D} \sum \limits_{i=1}^{N_D} \left [ \hat{u}(\mathbf{x}_i^D;\theta) - g^{*}_i \right ]^2 +\frac{1}{N_N}\sum \limits_{i=1}^{N_N} f_N(\mathbf{x}_i^N;\gamma,\theta)^2\\
    &+ \frac{1}{N_c} \sum \limits_{i=1}^{N_c} f(\mathbf{x}^c_i; \gamma, \theta) ^2.
  \end{split}
\end{equation}
The DNNs are trained by minimizing the loss function~(\ref{eq:loss_fn_linear}) as described in Section \ref{philms}. Throughout this work, we use feed-forward networks with two hidden layers and 50 units per layer.

\begin{figure}[htbp]
  \centering%
  \includegraphics[width=0.45\textwidth]{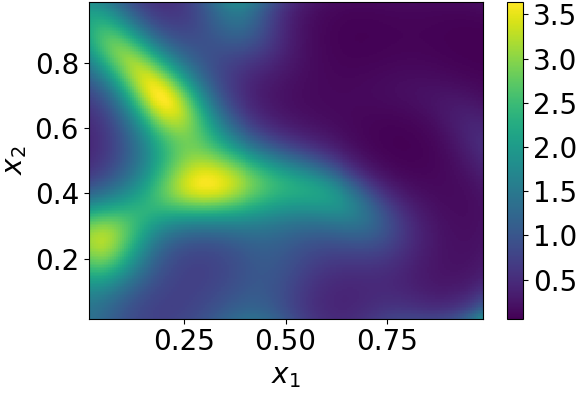}%
  \includegraphics[width=0.45\textwidth]{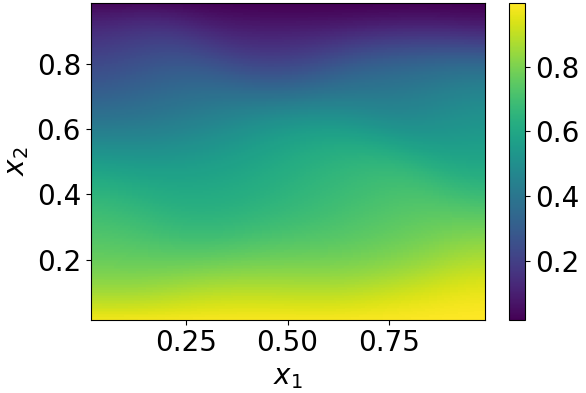}%
  \caption{Reference $K$ (left) and $u$ (right) fields.}
  \label{fig:references}
\end{figure}

\begin{figure}[htbp]
  \centering
   \begin{subfigure}[t]{0.49\textwidth}
  \includegraphics[width=0.99\textwidth]{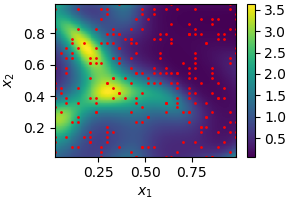}
    \caption{ $\hat{K}(\mathbf{x})$}  
 \end{subfigure}  
    \begin{subfigure}[t]{0.49\textwidth}
    \includegraphics[width=0.99\textwidth]{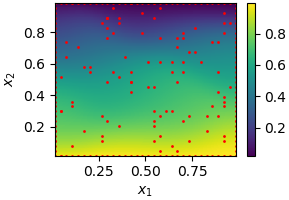}
      \caption{ $\hat{u}(\mathbf{x})$}  
 \end{subfigure}  
    \begin{subfigure}[t]{0.49\textwidth}
      \includegraphics[width=0.99\textwidth]{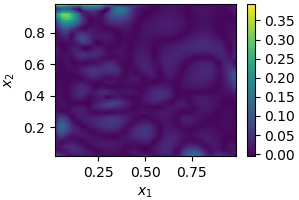}
        \caption{ $|{K}(\mathbf{x}) - \hat{K}(\mathbf{x})|$}  
 \end{subfigure}  
    \begin{subfigure}[t]{0.49\textwidth}
        \includegraphics[width=0.99\textwidth]{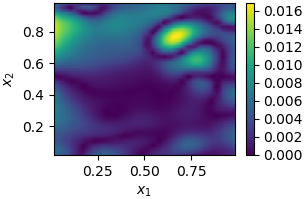}
         \caption{ $|{u}(\mathbf{x}) - \hat{u}(\mathbf{x})|$}  
 \end{subfigure}        
  \caption{Estimated (a) $\hat{K}(\mathbf{x})$ and (b) $\hat{u}(\mathbf{x})$ fields. The red dots indicate the observation locations. Absolute point errors in (c) $\hat{K}(\mathbf{x})$ and (d) $\hat{u}(\mathbf{x})$. $N_K = 250$, $N_u = 100$, and $N_c = 1024$.}
  \label{fig:est_with_BC}
\end{figure}

To demonstrate the proposed approach, we generate a reference $\ln K(\mathbf{x})$ field as a realization of the Gaussian process with zero mean and covariance function $C(\mathbf{x}, \mathbf{x}')=\sigma^2 \exp(-\| \mathbf{x}-\mathbf{x}' \|^2 / 2\lambda^2)$, with $\sigma = 1$ and $\lambda = 0.15$.
The reference $u$ is generated by solving Eqs (\ref{eq:linearPDE})--(\ref{eq:linearPDE-bc-N}) using the finite volume (FV) method with the two-point flux approximation and a cell-centered regular mesh with 1024 cells.
Figure~\ref{fig:references} presents the reference $K$ and $u$ fields.
We randomly choose $N_K$ and $N_u$ FV cell centroids as measurement locations for $K$ and $u$, respectively.
These measurement locations are shown in Figure~\ref{fig:est_with_BC}.
For evaluating the loss function and training the DNNs, we use $N_c = 1024$ uniformly distributed collocation points.

  We quantify the error between estimated and reference $K$ and $u$ fields in terms of the relative $L_2$ errors, defined as
  \begin{equation*}
    \varepsilon_u = \frac{\int_\Omega \left [ u(\mathbf{x})-\hat{u}(\mathbf{x}) \right ]^2 \mathrm{d} \mathbf{x}}{\int_\Omega u^2(\mathbf{x}) \, \mathrm{d} \mathbf{x}}, \quad \varepsilon_K = \frac{\int_\Omega \left [ K(\mathbf{x})-\hat{K}(\mathbf{x}) \right ]^2 \, \mathrm{d} \mathbf{x}}{\int_\Omega K^2(\mathbf{x}) \, \mathrm{d} \mathbf{x}}.
  \end{equation*}

Figure~\ref{fig:est_with_BC} shows the estimated $\hat{K}$ and $\hat{u}$ with $N_K = 250$, $N_u = 100$, and $N_c = 1024$. %
The relative $L_2$ errors %
are $\varepsilon_u \approx 0.5$\% and $\varepsilon_K \approx 1.7$ \%.
Figure~\ref{fig:est_with_BC} also depicts the point-wise absolute error in estimates of $K$ and $u$.
The point errors in $K$ are concentrated in the upper left corner where no $K$ measurements are available with a maximum point error of approximately $30\%$.
The point errors in $u$ are much smaller (maximum error is approximately $1\%$) and more uniformly distributed throughout the domain. 

\begin{figure}[htbp]
  \centering%
  \includegraphics[width=0.45\textwidth]{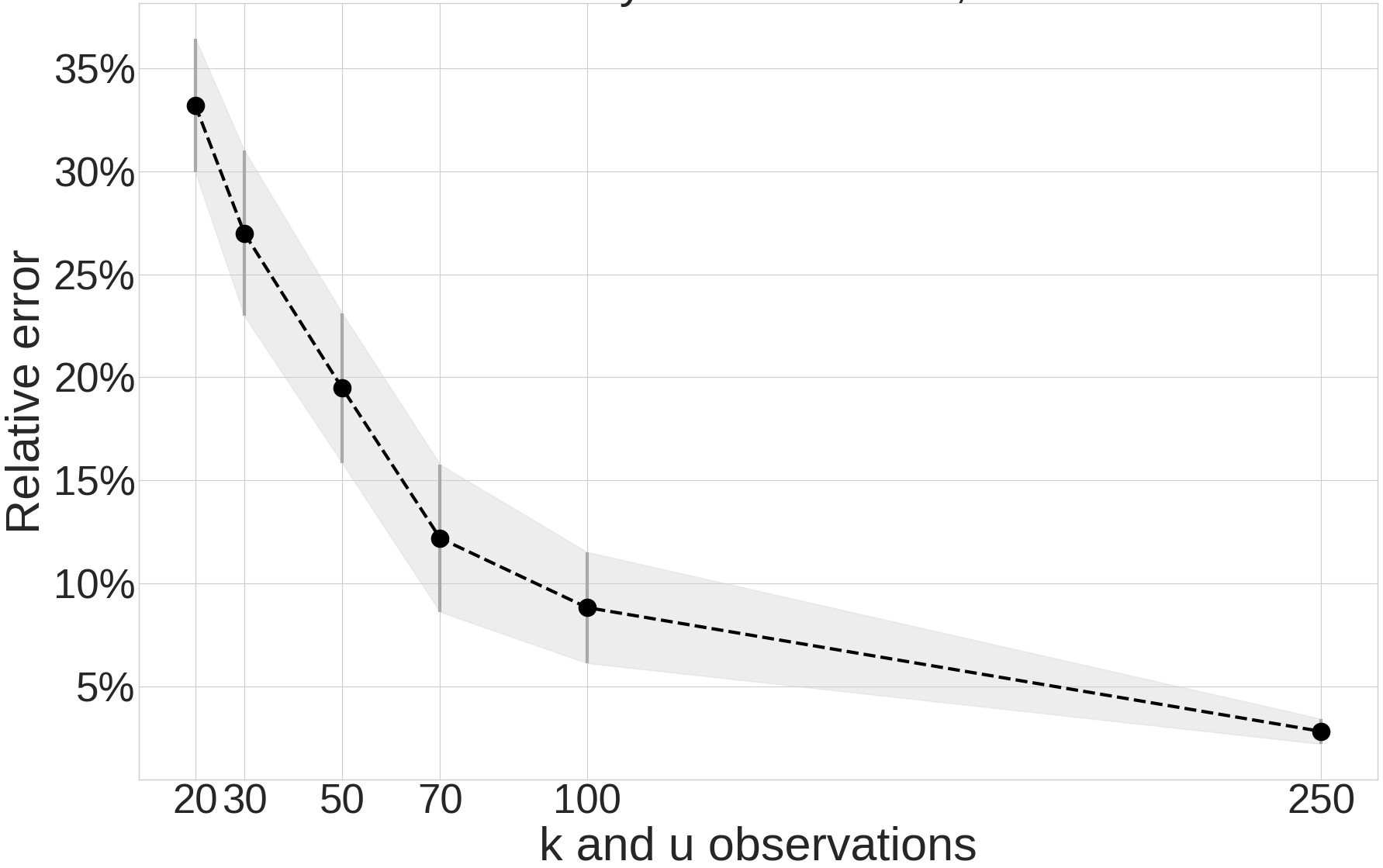}%
  \includegraphics[width=0.45\textwidth]{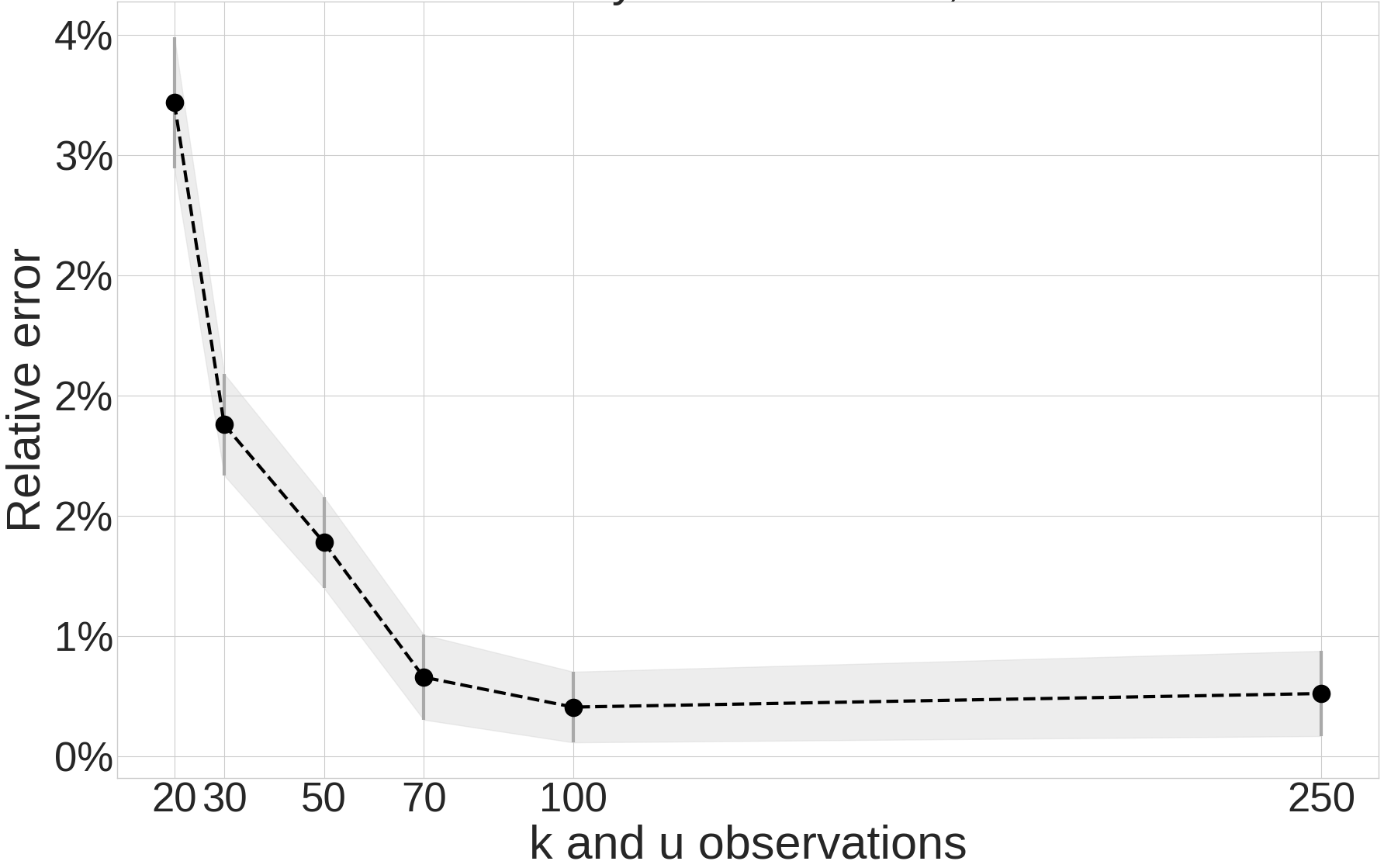}
  \caption{Mean and standard deviation of $\hat{K}$ and $\hat{u}$ obtained with 11 different network initializations using Xavier's initialization scheme as a function of $N=N_K=N_u$. The $u$ and $K$ measurement locations are fixed, and $N_c = 1024$ collocation points are used.}
  \label{fig:initialization}
\end{figure}

Next, we study the effect of DNN initialization on the estimated $K$ and $u$.
  For this purpose, we draw multiple initializations of the DNNs employing Xavier's scheme (see Section~\ref{philms}), and, for each initialization, we train the DNNs.
  We quantify the effect of initialization in terms of the mean and standard deviation of the relative $L_2$ errors of the estimated fields obtained for each initialization, i.e.,
  \begin{equation*}
    \overline{\varepsilon}_{(\cdot)}  = \frac1{N_s} \sum_{i=1}^{N_s} \varepsilon_{(\cdot), i}, \quad \sigma_{\varepsilon_{(\cdot)}} =  \sqrt{\frac1{N_s} \sum_{i=1}^{N_s} \left ( \varepsilon_{(\cdot), i} - \overline{\varepsilon}_{(\cdot)} \right )^2},
  \end{equation*}
  where $\varepsilon_{(\cdot), i}$ is the relative $L_2$ error for either $K$ or $u$ for the $i$th initialization, and $N_s$ is the number of network initializations.
Figure \ref{fig:initialization} shows the mean and standard deviation of $\varepsilon_u$ and $\varepsilon_K$
as a function of $N=N_K=N_u$ obtained from $N_s=11$ different network initializations.
 The $u$ and $K$ measurement locations are the same in these simulations, and $N_c = 1024$ collocation points are used.
As before, we see that uncertainty in $\hat{u}$ is much smaller than in $\hat{K}$, i.e., $\sigma_{\varepsilon_K} >> \sigma_{\varepsilon_u}$.
For both $\hat{u}$ and $\hat{K}$, the standard deviation associated with the initialization is approximately 10 times smaller than the mean value (the coefficient of variation is $\approx 0.1$), which indicates that the initialization of the DNNs does not have a significant effect on DNN predictions. 

\begin{figure}[htbp]
  \centering%
 \begin{subfigure}[t]{0.49\textwidth}
   \includegraphics[width=0.99\textwidth]{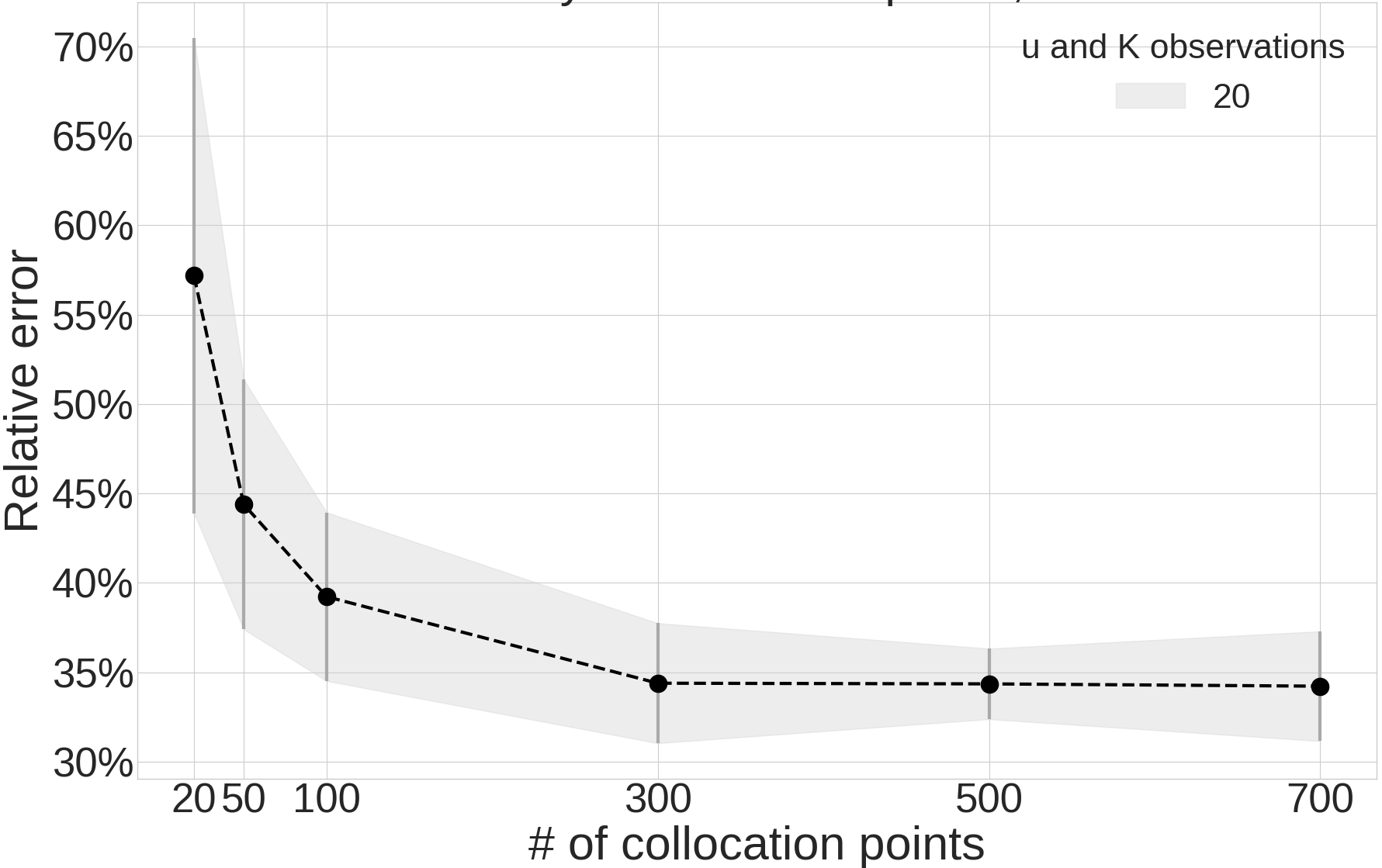}%
  \caption{ $\overline{\varepsilon}_{K}$ and $\sigma_{\varepsilon_{K}}$}  
 \end{subfigure}        
 \begin{subfigure}[t]{0.49\textwidth}  
  \includegraphics[width=0.99\textwidth]{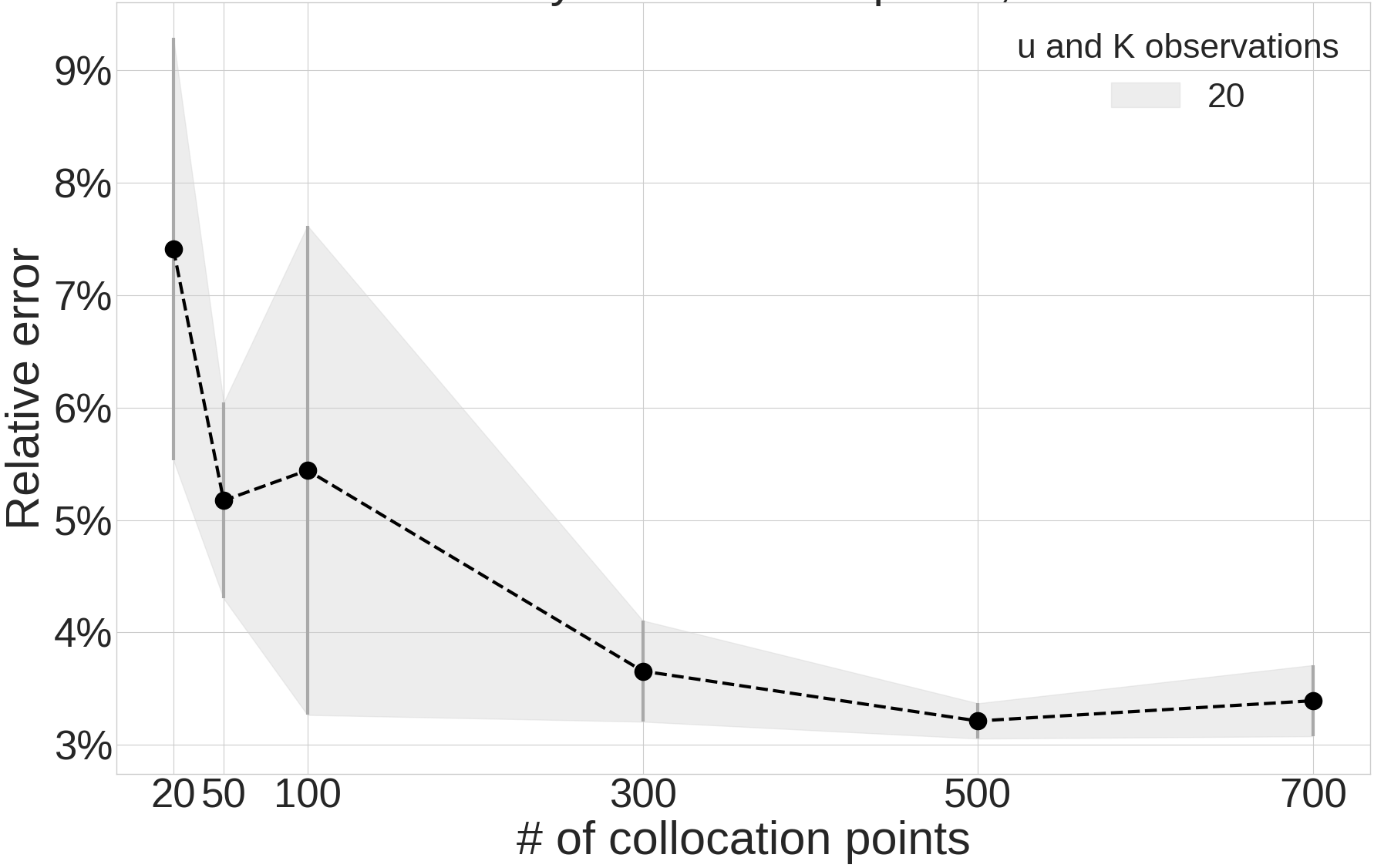}%\\
    \caption{$\overline{\varepsilon}_{u} $ and $\sigma_{\varepsilon_{u}}$}  
 \end{subfigure}        
 \begin{subfigure}[t]{0.49\textwidth}  
  \includegraphics[width=0.99\textwidth]{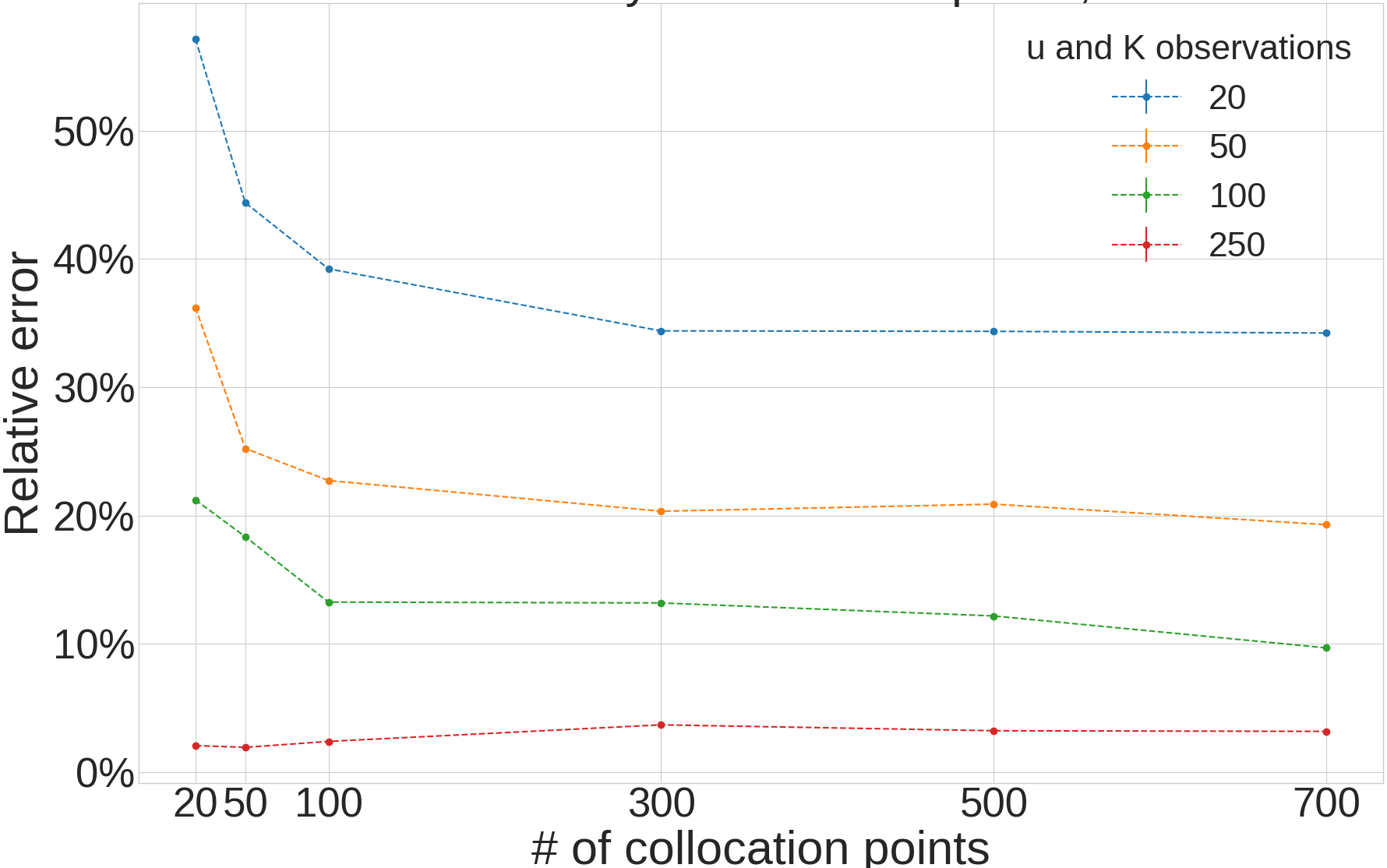}%
    \caption{ $\overline{\varepsilon}_{K}$}  
 \end{subfigure}        
 \begin{subfigure}[t]{0.49\textwidth}   
  \includegraphics[width=0.99\textwidth]{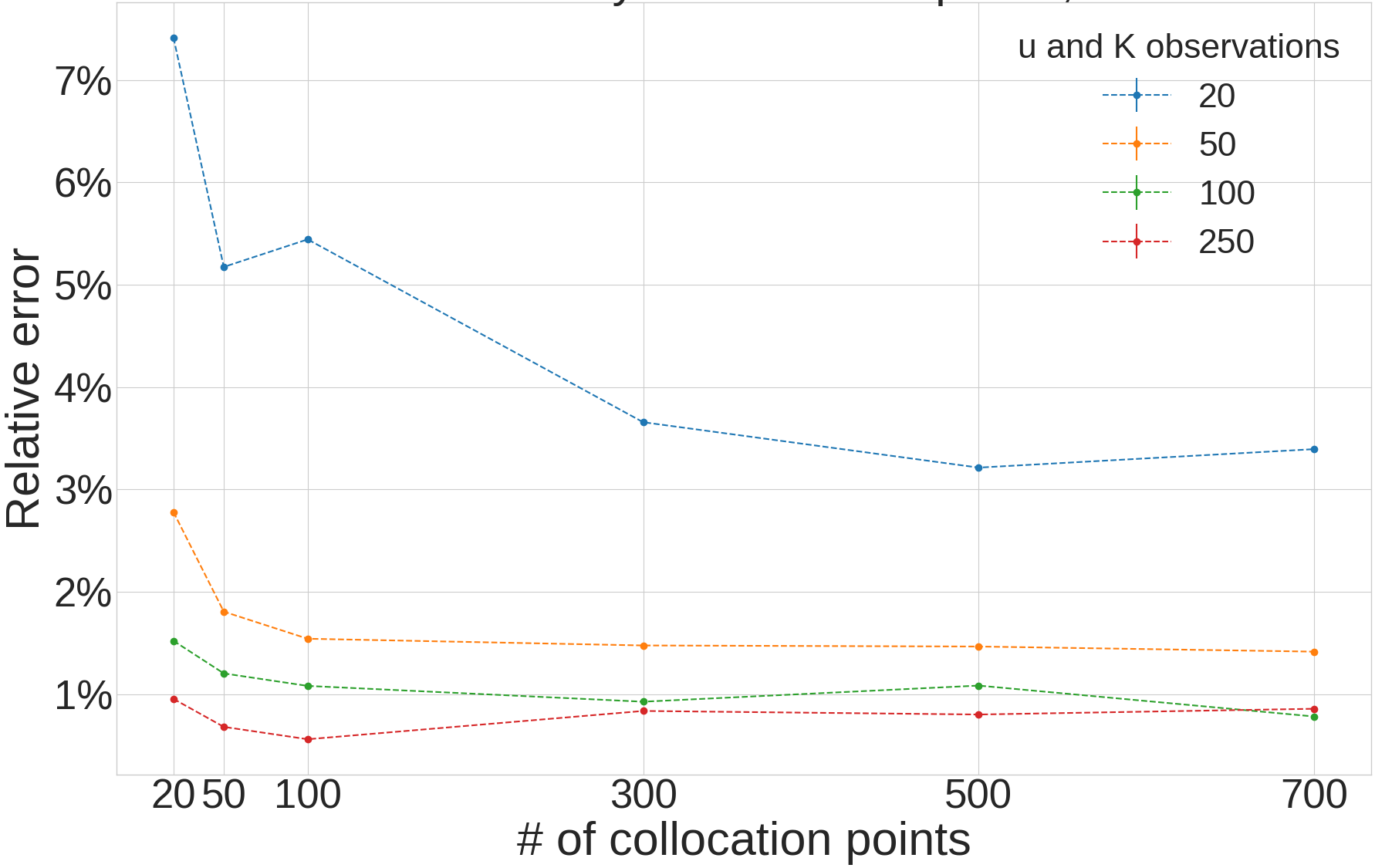}%
    \caption{ $\overline{\varepsilon}_{u}$}  
  \end{subfigure}       
  \caption{(a)-(b): Mean and standard deviation of the estimated $u$ and $K$ relative errors as a function of the number of collocation points $N_c$. For a given $N_c$, the DNNs are trained 11 times for different configurations of collocation points to compute the mean and variance of the relative errors. The $u$ and $K$ measurement locations are fixed, and $N=N_K=N_u=20$. Shaded area width is equal to two standard deviations.  (c)-(d): Mean $K$ and $u$ relative errors as a function of $N_c$ and $N$.}
  \label{fig:collocation}
\end{figure}

For the DNNs training, the governing PDEs are enforced at $N_c$ collocation points. To study the effect of the number and location of collocation points, we compute the relative errors $\varepsilon_u$ and $\varepsilon_K$ as a function of the number and location of the collocation points. Figures \ref{fig:collocation}(a) and (b) show the mean and standard deviation of the relative errors versus $N_c$ for $N=N_u=N_K=20$.
For a given $N_c$, the DNNs are trained $N_s=11$ times for different locations chosen via 
Latin hypercube sampling \cite{mckay1979comparison} to compute the mean and variance of the relative errors.
The error in $\hat{K}$ is about eight times larger than in $\hat{u}$. 
As expected, the mean and standard deviation of $\varepsilon_u$ and $\varepsilon_K$ decrease with increasing $N_c$ until they reach asymptotic values at approximately $N_c = 300$, which is approximately 33\% of the number of grid points in these simulations.
The location of collocation points has a notable  effect on the errors, especially for a relatively small $N_c$, which is evident from relatively large coefficients of variation  
$\sigma_{\varepsilon_{K}} / \overline{\varepsilon}_{K}$ 
and
$\sigma_{\varepsilon_{u}} / \overline{\varepsilon}_{u}$. 

Figures \ref{fig:collocation}(c) and (d) show that  $\overline{\varepsilon}_u$ and $\overline{\varepsilon}_K$ asymptotically decrease for four considered values of $N$. The asymptotic values of $\overline{\varepsilon}_u$ and $\overline{\varepsilon}_K$ also decrease with increasing $N$. 
For all considered $N$, imposing PDE constraints reduces the mean error in $K$ and $u$ by close to 50\%.
Here, as in Figures \ref{fig:collocation}(a) and (b), $\varepsilon_u$ is significantly smaller than $\varepsilon_K$.

\begin{figure}[htbp]
  \centering%
   \begin{subfigure}[t]{0.49\textwidth}
  \includegraphics[width=0.99\textwidth]{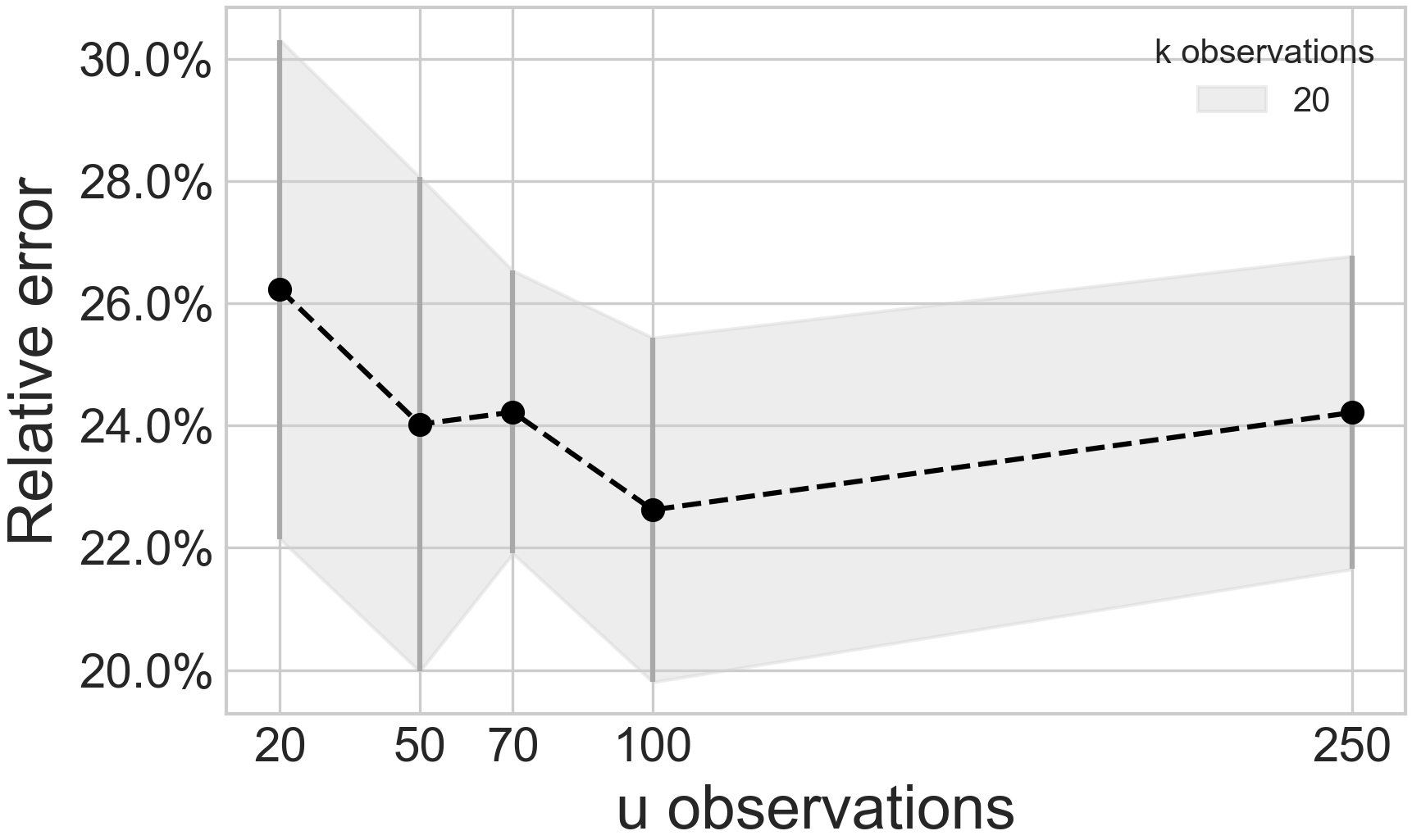}%
   \caption{$\overline{\varepsilon}_{K}$ and $\sigma_{\varepsilon_{K}}$}  
 \end{subfigure}        
 \begin{subfigure}[t]{0.49\textwidth}   
  \includegraphics[width=0.99\textwidth]{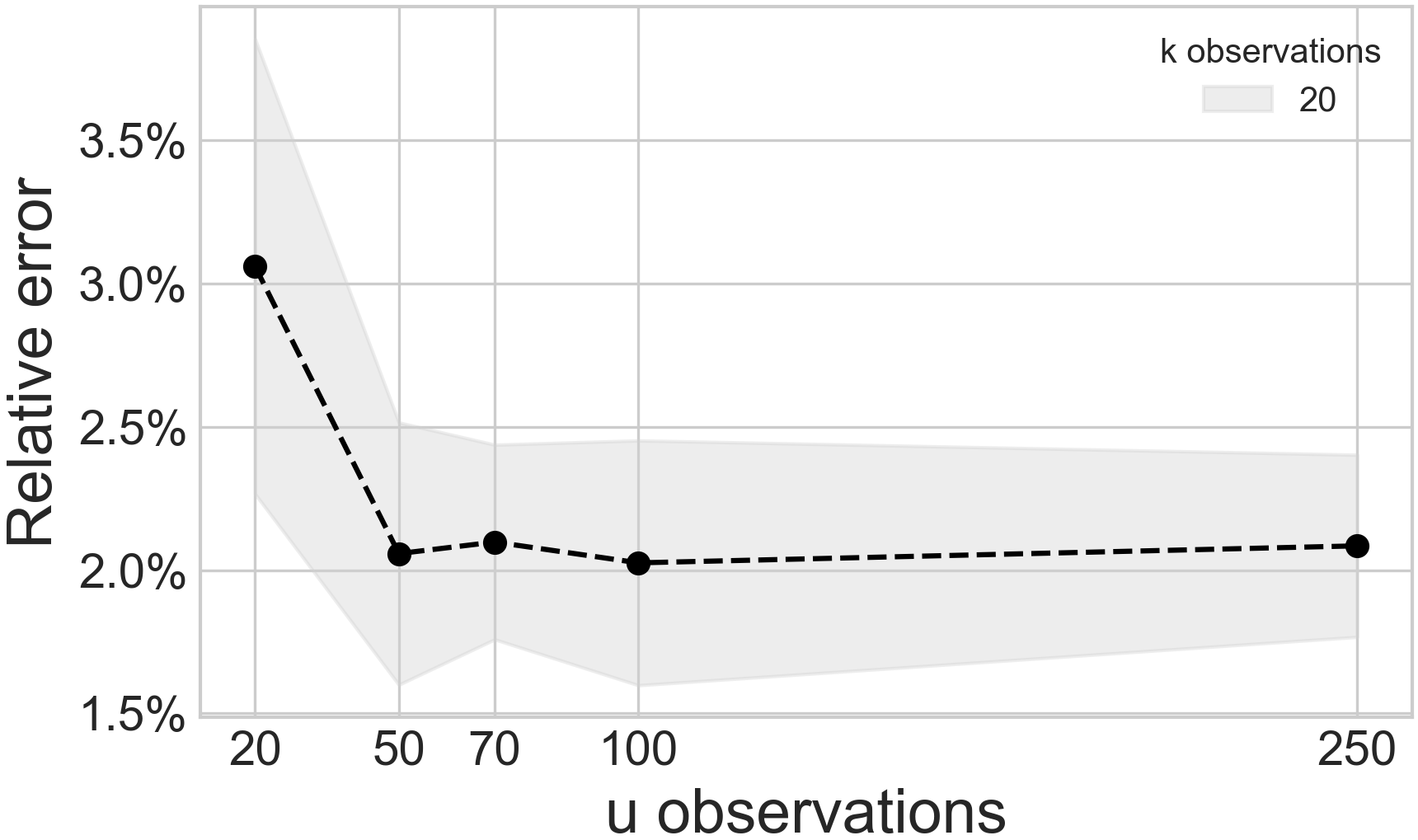}
   \caption{$\overline{\varepsilon}_{u}$ and $\sigma_{\varepsilon_{u}}$}  
 \end{subfigure}        
 \begin{subfigure}[t]{0.49\textwidth}   
  \includegraphics[width=0.99\textwidth]{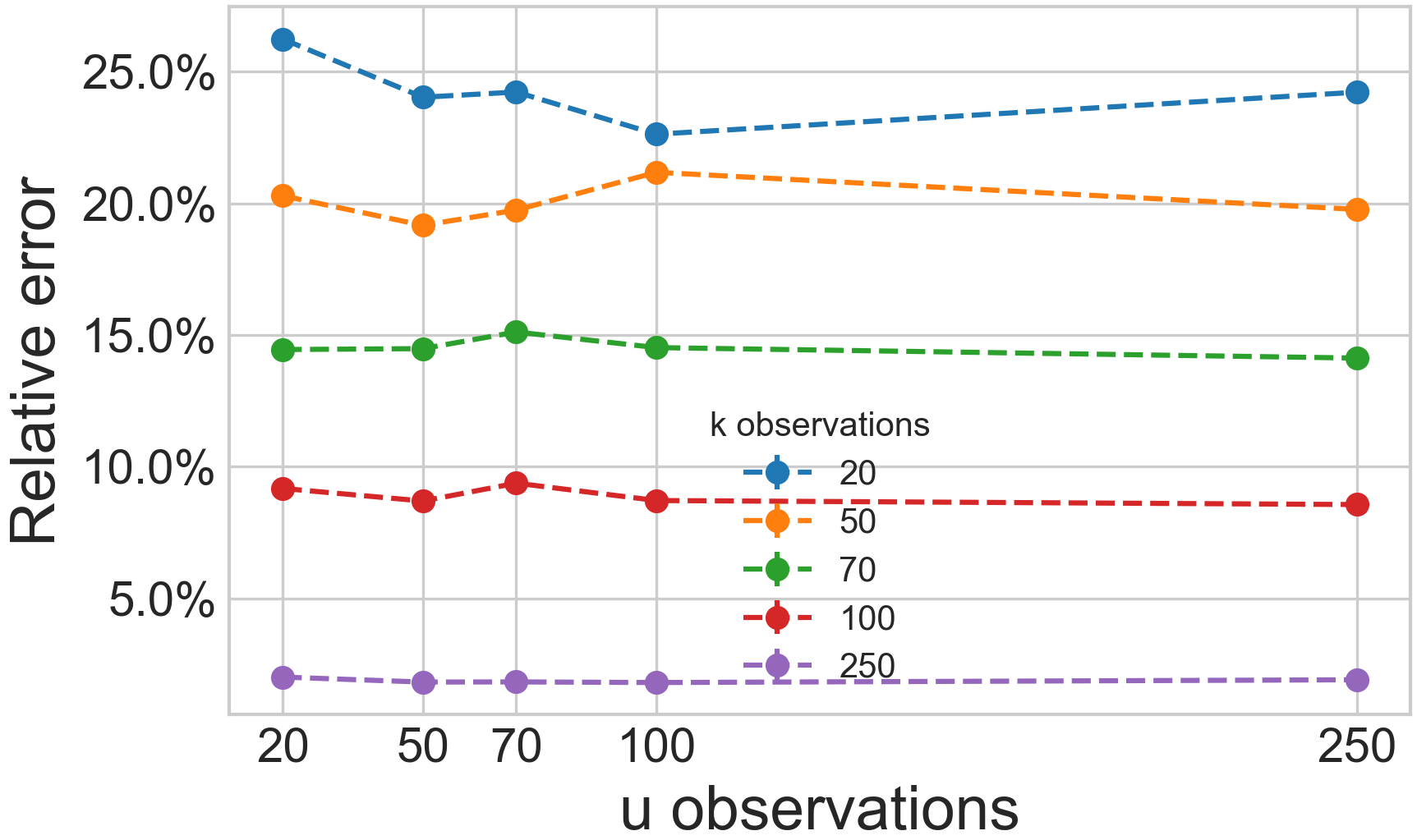}%
    \caption{$\overline{\varepsilon}_{K}$}  
 \end{subfigure}        
 \begin{subfigure}[t]{0.49\textwidth}  
  \includegraphics[width=0.99\textwidth]{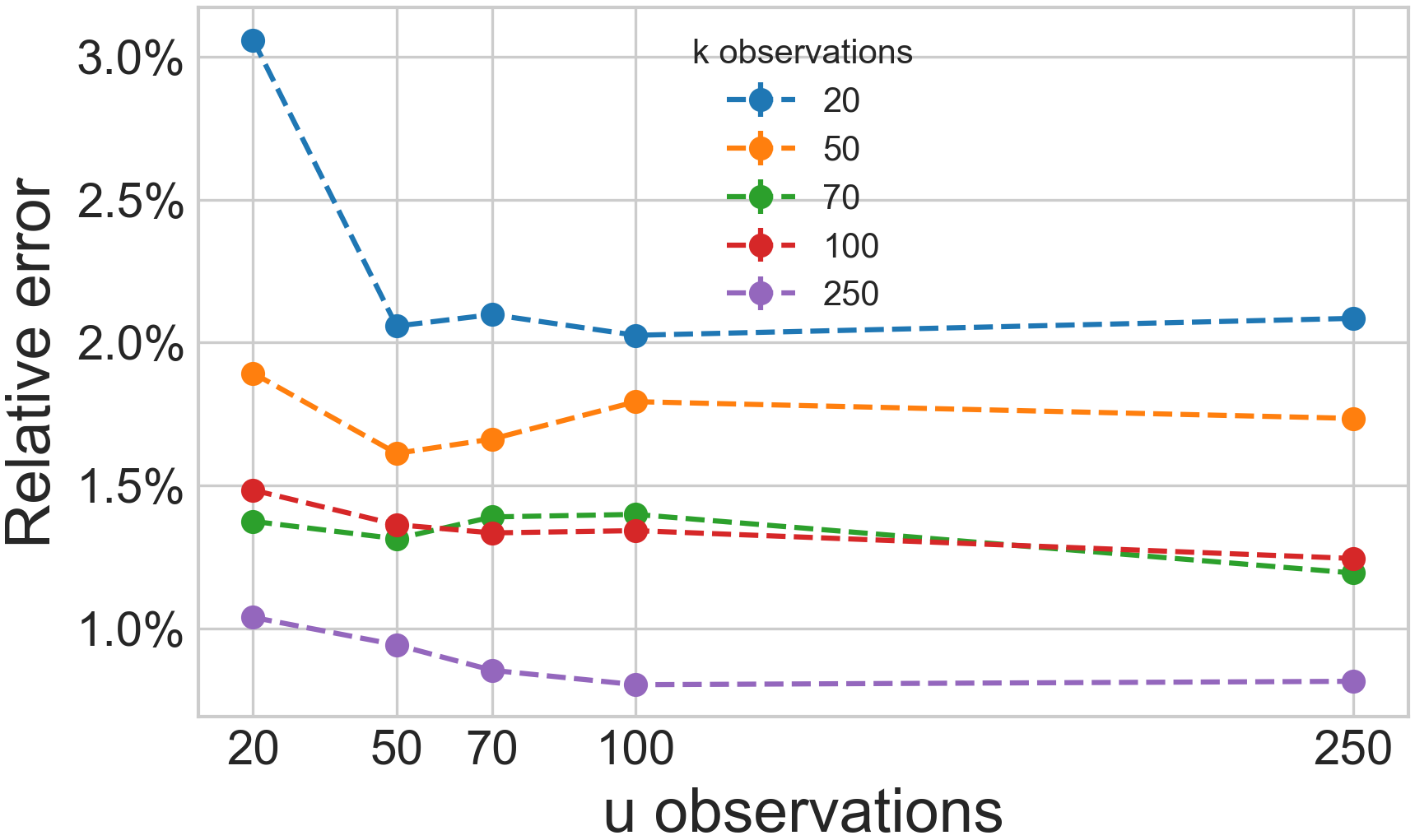}%
    \caption{$\overline{\varepsilon}_{u}$}  
 \end{subfigure}     
  \caption{ (a) and (b): Mean relative $L_2$ error of the predicted $K$ and  $u$ as a function of the number of $u$ observations. The number of $K$ observations is 20. Shaded area width is equal to two standard deviations computed for 11 different configurations of $u$ observations. (c) and (d):  Mean relative $L_2$ error of the predicted $K$ and  $u$ as a function of the number of $u$ and $K$ observations.}
  \label{fig:fixed_k}
\end{figure}

\begin{figure}[htbp]
  \centering
     \begin{subfigure}[t]{0.49\textwidth}
 \includegraphics[width=0.99\textwidth]{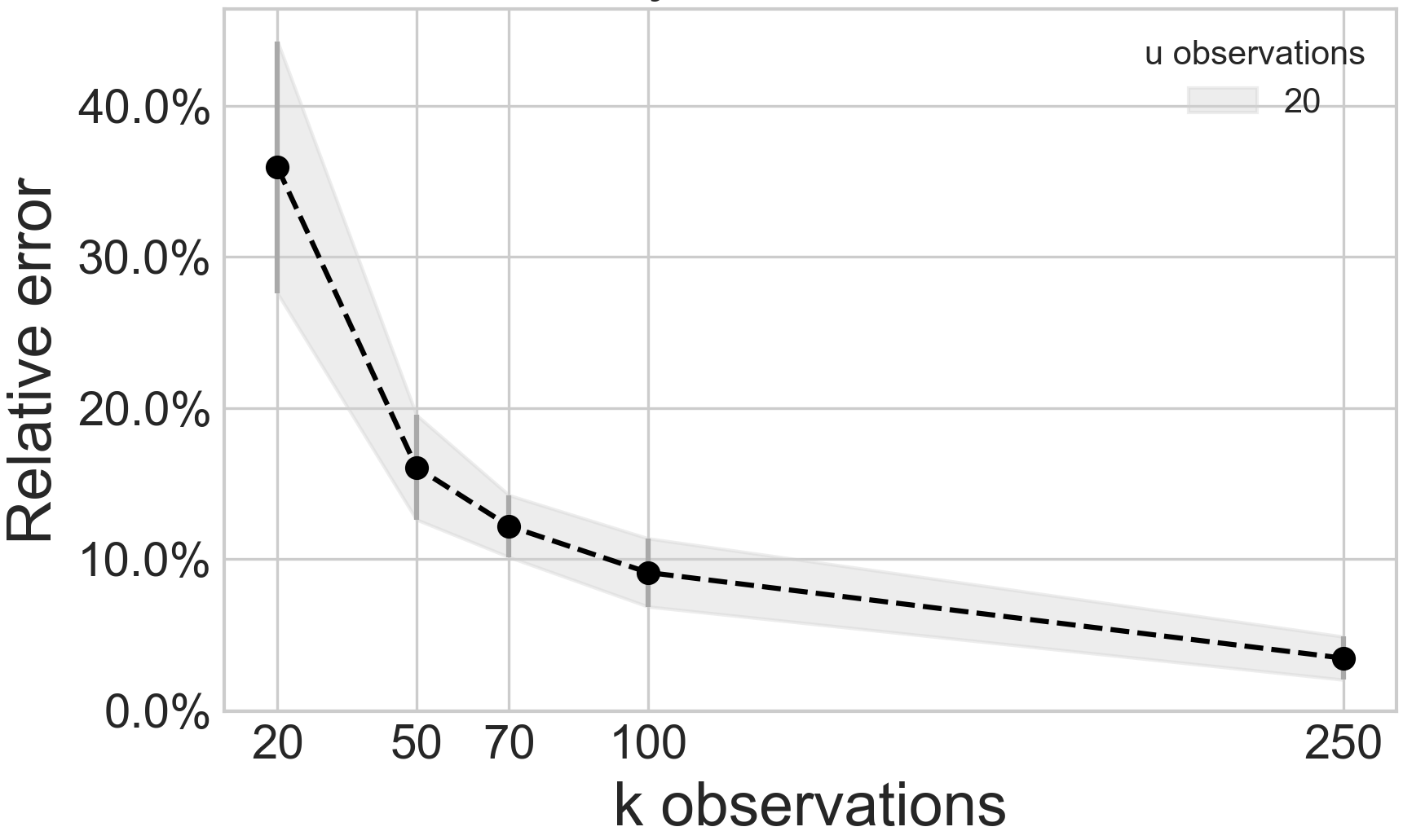}
    \caption{$\overline{\varepsilon}_{K}$ and $\sigma_{\varepsilon_{K}}$}  
 \end{subfigure}        
 \begin{subfigure}[t]{0.49\textwidth}  
 \includegraphics[width=0.99\textwidth]{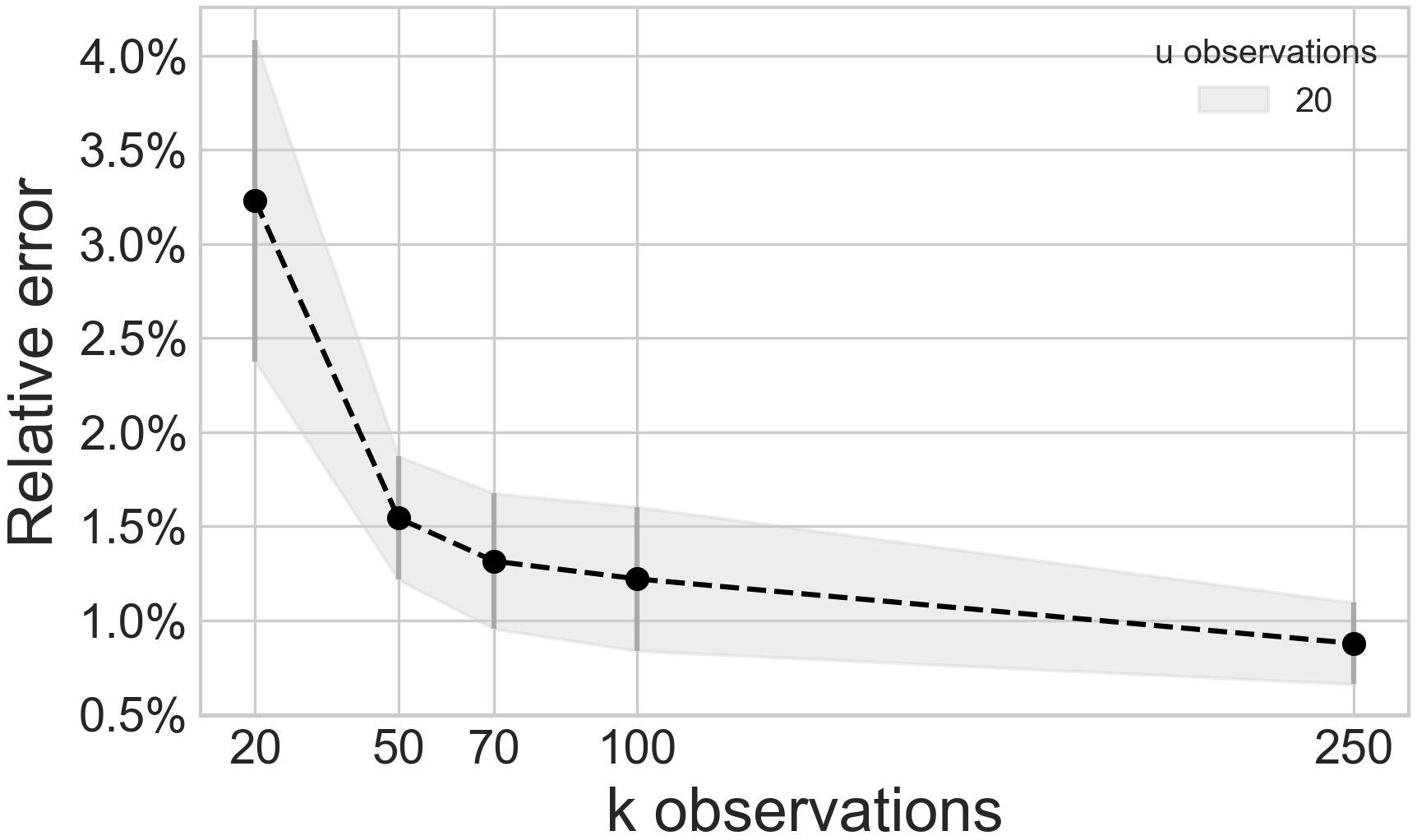}
    \caption{$\overline{\varepsilon}_{u}$ and $\sigma_{\varepsilon_{u}}$}  
 \end{subfigure}        
 \begin{subfigure}[t]{0.49\textwidth}  
 \includegraphics[width=0.99\textwidth]{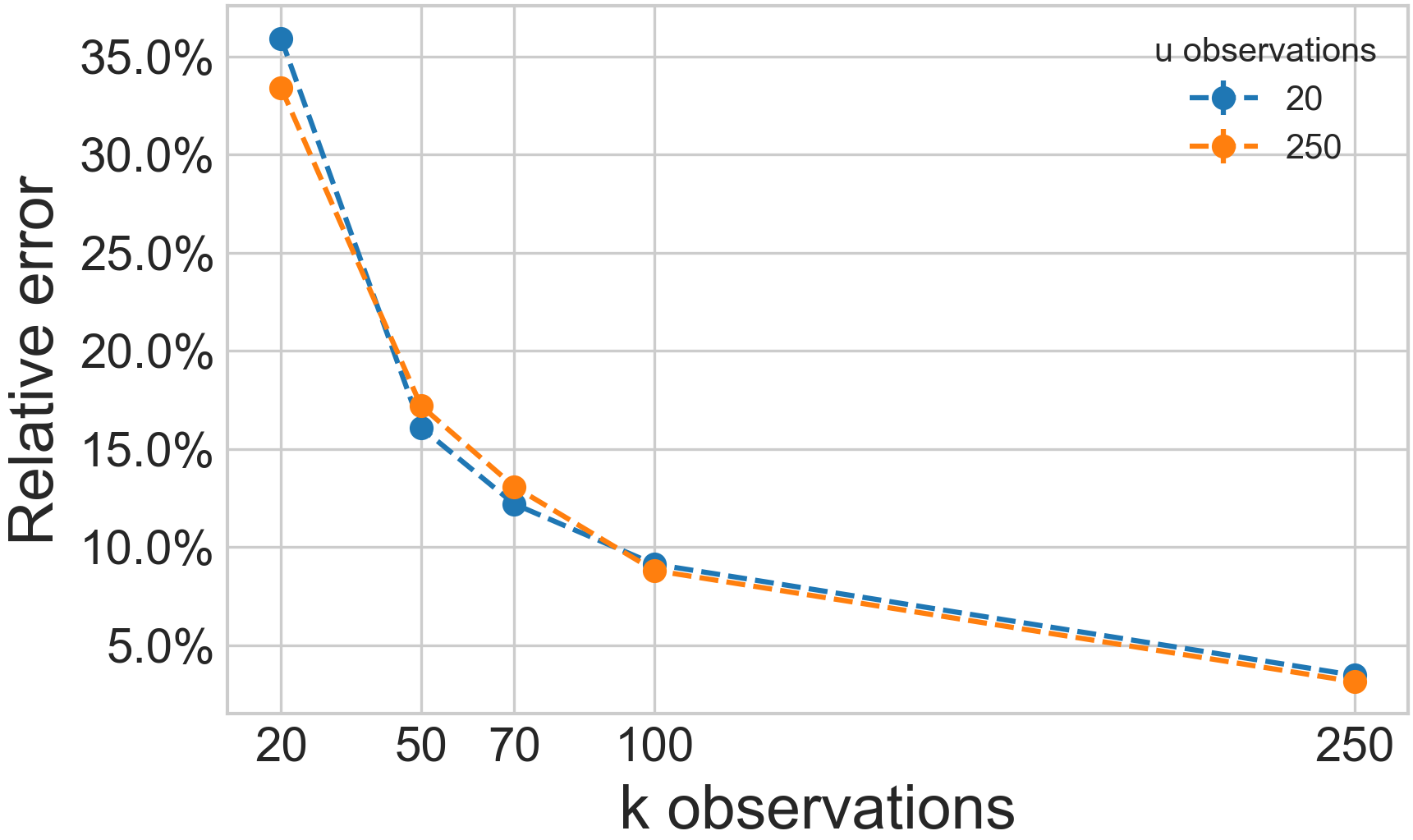}
    \caption{$\overline{\varepsilon}_{K}$}  
 \end{subfigure}        
 \begin{subfigure}[t]{0.49\textwidth}  
 \includegraphics[width=0.99\textwidth]{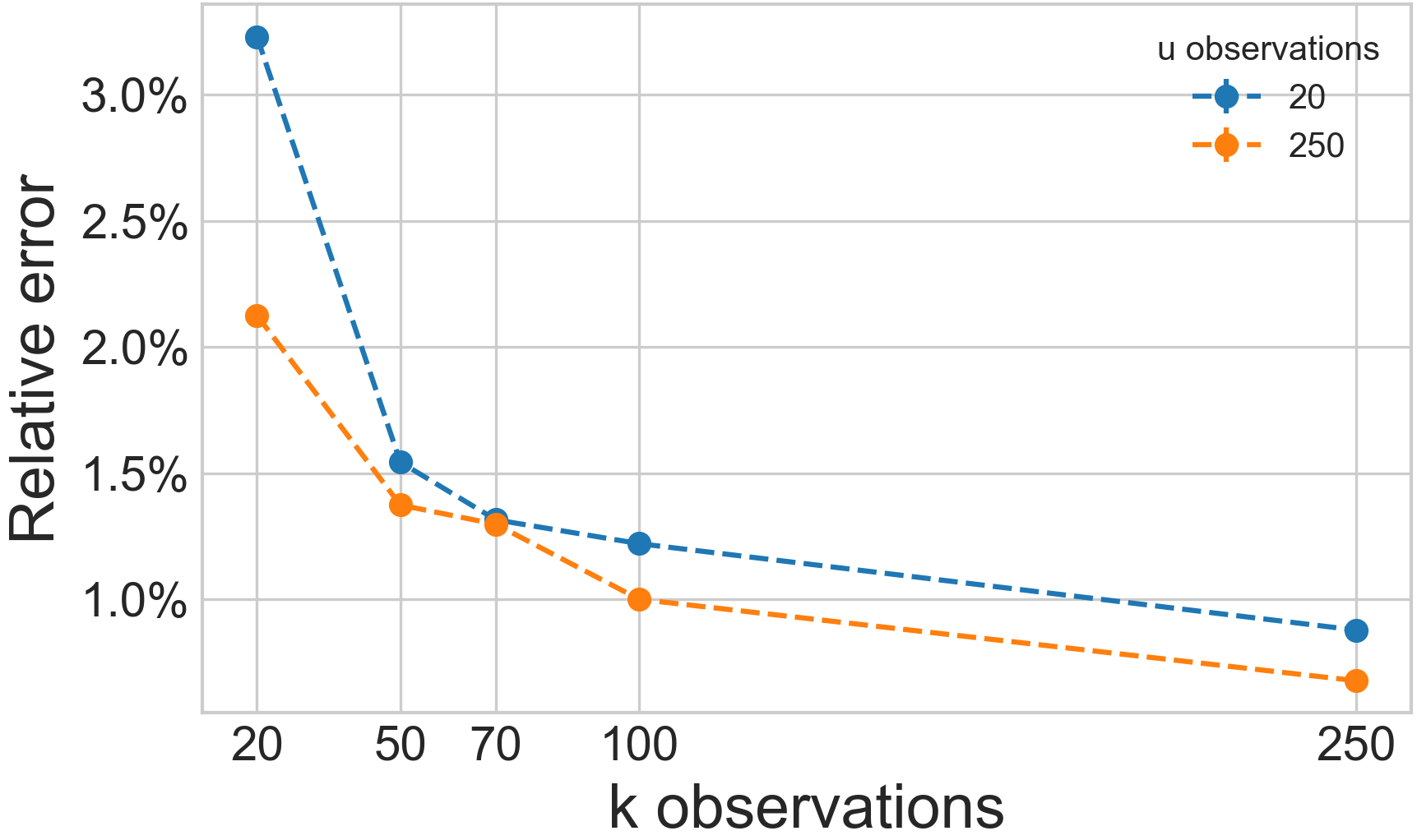}
     \caption{$\overline{\varepsilon}_{u}$}  
 \end{subfigure}    
 \caption{(a) and (b): Mean relative $L_2$ error of the predicted $K$ and  $u$ as a function of the number of $K$ observations. The number of $u$ observations is 20. Shaded area width is equal to two standard deviations computed for 11 different configurations of $K$ observations. (c) and (d):  Mean relative $L_2$ error of the predicted $K$ and  $u$ as a function of the number of $u$ and $K$ observations.}
  \label{fig:fixed_u}
\end{figure}

In Figures \ref{fig:fixed_k} and \ref{fig:fixed_u}, we study how the number of $K$ observations versus the number of $u$ observations affects $\varepsilon_u$ and $\varepsilon_K$. Here, the number of collocation points is $N_c = 1024$. Figures \ref{fig:fixed_k}(a) and (b) depict the effect of $N_u$ on the mean and standard deviation of $\varepsilon_u$ and $\varepsilon_K$ for $N_K = 20$. For each $N_u$ value, we estimate $\hat{K}$ and $\hat{u}$ with $N_s=11$ different distributions of the $u$ measurement  locations generated with the Latin hypercube sampler. Then, we compute $\overline{\varepsilon}_u$, $\overline{\varepsilon}_K$, $\sigma_{\varepsilon_u}$, and $\sigma_{\varepsilon_K}$. We see that $N_u$ does not significantly affect $\varepsilon_K$ with $\overline{\varepsilon}_K \approx24$\% and $ \sigma_{\varepsilon_K} \approx6$\%. For $\hat{u}$, $\overline{\varepsilon}_u$ deceases from more than 3\% for $N_u=20$ to $\approx 2$\% for $N_u>100$. The standard deviation $ \sigma_{\varepsilon_u}$ is
1.5\% for $N_u=20$ and less than 1\% for $N_u>50$.

Figures \ref{fig:fixed_k} (c) and (d) show $\overline{\varepsilon}_K$ and $\overline{\varepsilon}_u$ as a function of $N_u$ for different $N_K$.
For all considered $N_K$, $\overline{\varepsilon}_K$ is practically independent of $N_u$ and decreases with increasing $N_K$. On the other hand, $\overline{\varepsilon}_u$ decreases with increasing $N_u$ and/or $N_K$.   

Figures \ref{fig:fixed_u}(a) and (b) reveal that all  $\overline{\varepsilon}_u$, $\overline{\varepsilon}_K$, $\sigma_{\varepsilon_u}$, and $\sigma_{\varepsilon_K}$ decrease with increasing $N_K$ for fixed $N_u$. Alternatively, Figures \ref{fig:fixed_u}(c) and (d) demonstrate that $N_u$ has a relatively minor effect on  $\overline{\varepsilon}_u$ and $\overline{\varepsilon}_K$. In all considered cases, $\overline{\varepsilon}_u$ is almost an order of magnitude smaller than $\overline{\varepsilon}_K$. 

The main conclusion to be drawn from  Figures \ref{fig:fixed_k} and \ref{fig:fixed_u} is that $K$ measurements are more important than $u$ measurements for reducing error in $\hat{K}$ and $\hat{u}$. This may be explained partially by the fact that $u(\mathbf{x})$ is much smoother than $K(\mathbf{x})$, and a relatively small number of $u$ measurements are needed to describe the $u$ field. Beyond this number (approximately 50 for this example), additional $u$ measurements do not have a significant effect on $\overline{\varepsilon}_u$ and $\overline{\varepsilon}_K$.

\begin{figure}[htbp]
  \centering
  \includegraphics[width=0.75\textwidth]{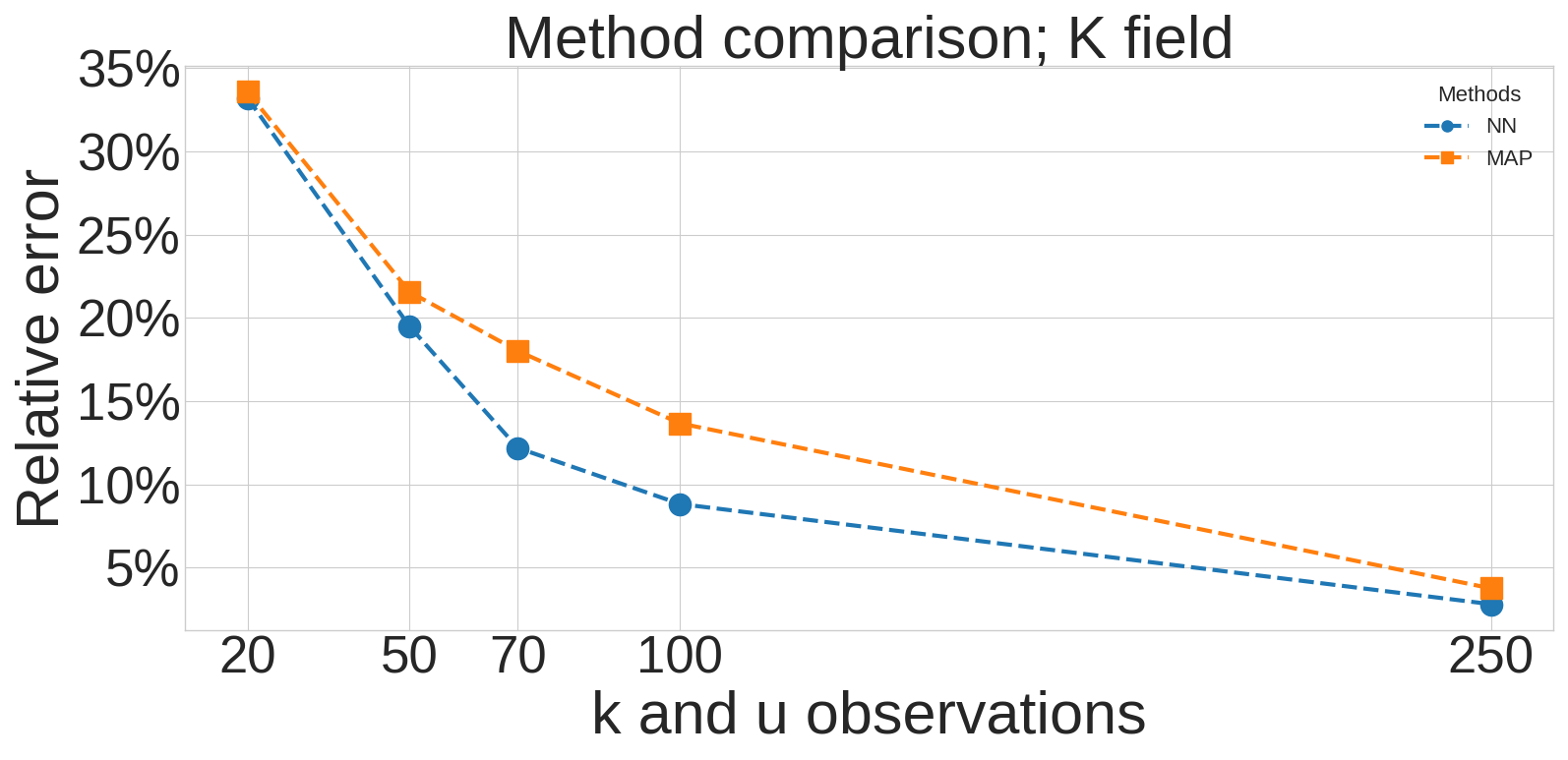}
  \caption{Comparison between $L_2$ error in the MAP estimate of $K$ and mean $L_2$ error $\overline{\varepsilon}_K$  in the DNN $K$ estimate.}
  \label{fig:mapvsnn}
\end{figure}

Finally, we compare our approach with the MAP method.
To regularize the $K$ estimate, we penalize the norm of the $K$ gradient, which prefers a smoother estimate \cite{barajassolano-2014-linear}.
For this regularizer and the regular FV discretization with $M = 1024$ cells, the MAP estimate $\hat{\mathbf{k}} \in \mathbb{R}^M$ is defined as the vector of cell-centered values of $K$ computed as the solution to the minimization problem
\begin{equation}
  \begin{aligned}
    & \hat{\mathbf{k}} = \argmin_{\mathbf{k}} & & \| \mathbf{u}^{*} - \mathbf{H}_u \mathbf{u} ||^2_2 + ||\ln \mathbf{k}^{*} - \mathbf{H}_K \ln \mathbf{k} ||_2^2+ \gamma \| \mathbf{L} \ln \mathbf{k} \|^2_2\\
    & \text{subject to } & & \mathbf{l}(\mathbf{u}, \mathbf{k}) = \bm{0},
  \end{aligned}
  \label{MAP}
\end{equation}
where $\mathbf{k}^{*} \equiv (K^{*}_1, \dots, K^{*}_{N_K})^{\top}$ and $\mathbf{u}^{*} \equiv (u^{*}_1, \dots, u^{*}_{N_u})^{\top}$ are the vectors of $K$ and $u$ observations, respectively; $\mathbf{H}_K \in \mathbb{R}^{N_K \times M}$ and $\mathbf{H}_u \in \mathbb{R}^{N_u \times M}$ are observation operators; $\mathbf{u}$ is the discretized state; $\mathbf{l}(\mathbf{u}, \mathbf{k}) = \bm{0}$ is the discretized problem~(\ref{eq:linearPDE})--(\ref{eq:linearPDE-bc-N}); $\mathbf{L}$ is the discrete gradient operator; and $\gamma<1$ is a problem-dependent coefficient.
The minimization problem (\ref{MAP}) is solved via the Levenberg-Marquardt (LM) algorithm~\cite{barajassolano-2014-linear}.

\begin{figure}[htbp]
  \centering
 \includegraphics[width=0.30\textwidth]{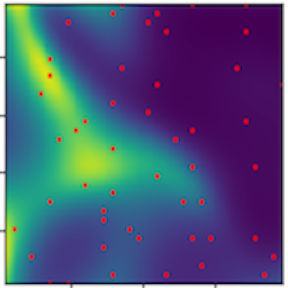}
 \includegraphics[width=0.30\textwidth]{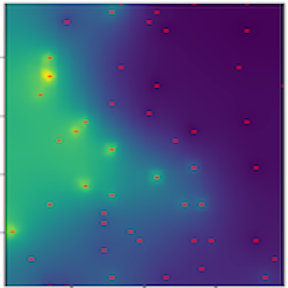}
 \caption{Comparison between the DNN prediction of $K$ (left) and the MAP prediction of $K$ (right) using 50 observations of both $u$ and $K$.}
  \label{fig:visualmap}
\end{figure}

For the considered problems, we find that the smallest $L_2$ error is obtained with $\gamma=10^{-6}$.
Figures \ref{fig:mapvsnn} and \ref{fig:visualmap} compare the two methods. Figure \ref{fig:mapvsnn} presents $\varepsilon_K$ as a function of $N=N_K=N_u$ for $\hat{K}$ found from the two methods. The physics informed DNNs produce $\hat{K}$ with smaller errors for all considered $N$. 

Figure \ref{fig:visualmap} depicts $\hat{K}$ obtained from the two methods for $N=50$. The $\hat{K}$ field estimated from the physics informed DNNs is significantly smoother (and closer to the reference $K$ shown in Figure \ref{fig:references}) than one estimated from MAP,  even though the $\varepsilon_K$ errors for the two methods are relatively  similar: $19\%$ for NN versus $22\%$ for MAP. The tent-like character  of the MAP prediction in Figure \ref{fig:visualmap} stems from the discrepancy with respect to observations is penalized more than the smoothness of the estimate (result of a relatively small $\gamma$). We find that larger $\gamma$ results in a smoother field, but also in a larger prediction error.

Both the MAP and physics-informed DNN methods involve an objective-function minimization.
In addition, MAP estimation via gradient-based optimization algorithms (such as the LM algorithm) requires computing the gradient of the predicted observations, $\mathbf{H}_u \mathbf{u}$ and $\mathbf{H}_K \ln \mathbf{k}$ with respect to $\mathbf{k}$.
For an FV discretization, this is done via the discrete adjoint method.
The total cost for each iteration of the LM algorithm is one forward solution of the PDE problem to evaluate the objective function and one adjoint solution to compute the gradient.
Therefore, MAP requires careful discretization of the PDE problem and formulation and solution of the corresponding adjoint problem.
In contrast, in the physics-informed DNN method, both spatial derivatives and the gradients with respect to DNN parameters are computed via automatic differentiation, and the methodology does not require solving the PDE problem or formulating an adjoint problem.

Finally, significant gains can be achieved in the physics informed DNN method performance by employing GPU accelerators.

\section{Nonlinear diffusion equation}\label{nonlinearPDE}

In this section, we consider a nonlinear diffusion equation with unknown state-dependent diffusion coefficient $K(u)$,
\begin{equation}
  \label{eq:nonlinearPDE}
  \nabla \cdot [K(u)\nabla u(\mathbf{x})] =0, \quad (x_1,x_2) \in (0,L_1) \times (0,L_2) 
\end{equation}
subject to the boundary conditions
\begin{align}
  \label{eq:nonlinearPDE-bc-D}
  u(\mathbf{x}) &= u_0, &&  x_1=L_1, \\
  \label{eq:nonlinearPDE-bc-N-1}
  -K(u) \frac{\partial u(\mathbf{x})}{\partial x_1} &= q, && x_1 = 0\\
  \label{eq:nonlinearPDE-bc-N-2}
  \frac{\partial u(\mathbf{x})}{\partial x_2} &= 0, && x_2 = \{ 0, L_2 \}.
\end{align}

This equation describes a two-dimensional horizontal unsaturated flow (flow of water and air) in a homogeneous porous medium, where $u(\mathbf{x})$ is the water pressure and $K(u)$ is the pressure-dependent partial conductivity of the porous medium~\cite{bear2013dynamics}.
In practice, $K(u)$ is difficult to measure directly.
Therefore, this work assumes that no measurements of $K(u)$ are available, and only $N_u$ measurements of $u$ are given.  

We define two DNNs for unknown $K(u)$ and $u(\mathbf{x})$,
\begin{equation}
  \hat{u}(\mathbf{x};\theta) = \mathcal{NN}_u(\mathbf{x}; \theta), \quad \hat{K}(u; \gamma) = \mathcal{NN}_K(u; \phi) 
\end{equation}
  and two auxiliary DNNs obtained by substituting the DNNs for $K$ and $u$ into (\ref{eq:nonlinearPDE}), (\ref{eq:nonlinearPDE-bc-N-1}), and~(\ref{eq:nonlinearPDE-bc-N-2}), 
\begin{gather*}
  f(\mathbf{x}; \theta, \gamma) = \nabla \cdot \left [ \mathcal{NN}_K \left ( \mathcal{NN}_u(\mathbf{x}; \theta); \gamma \right ) \nabla  \mathcal{NN}_u(\mathbf{x}; \theta)  \right ] = \mathcal{NN}_f (\mathbf{x}; \theta, \gamma), \\
  \mathbf{f}_N(\mathbf{x}; \phi, \gamma) = -\mathcal{NN}_K \left ( \mathcal{NN}_u(\mathbf{x}; \theta); \gamma \right ) \nabla \mathcal{NN}_u(\mathbf{x}; \theta) = \mathcal{N} \mathcal{N}_N (\mathbf{x}; \theta, \gamma),
\end{gather*}
  where $\mathbf{f}_N$ is a vector DNN with two DNN components $f^{(x_1)}_N$ and $f^{(x_2)}_N$.
Then, the loss function becomes
\begin{equation*}
  \begin{split}
    \mathcal{L}(\theta, \gamma) &=  \frac{1}{N_u} \sum \limits_{i=1}^{N_u} \left [ \hat{u}(\mathbf{x}^u_i; \theta) - u^*_i \right ]^2 \\
    &+ \frac{1}{N_c} \sum \limits_{i=1}^{N_c} f(\mathbf{x}_i; \theta, \gamma)^2 + \frac{1}{N_D} \sum \limits_{i=1}^{N_D} \left [ \hat{u}(\mathbf{x}_i^D; \theta) - u_0 \right ]^2\\
    &+ \frac{1}{N_N^{(x_1)}} \sum \limits_{i=1}^{N_N} \left [ f_N^{(x_1)}(\mathbf{x}_i^{N, x_1}; \phi, \gamma) - q \right ]^2 + \frac{1}{N_N^{(x_2)}} \sum \limits_{i=1}^{N_N} \left [f_N^{(x_2)}(\mathbf{x}_i^{N,x_2}; \phi, \gamma) \right]^2,
  \end{split}
  \label{loss_fn_Richards} 
\end{equation*}
where 
$\mathbf{x}_i^{N,x_1}$ $(i=1,...,N_N^{(x_1)})$ are the collocation points on the Neumann boundary $(x_1=0, x_2)$ and $\mathbf{x}_i^{N,x_2}$ $(i=1,...,N_N^{(x_2)})$ are the collocation points on the Neumann boundaries $(x_1, x_2=0)$ and $(x_1, x_2=L_2)$.

This model is tested with data generated using the Subsurface Transport Over Multiple Phases (STOMP) code \cite{white1995modeling} with the van Genuchten model \cite{van1980closed} for the $K(u)$ function
\begin{gather}
  \label{vGK}
  K(s(u)) = K_s s^{\frac12} \left[1- \left(1-s^{\frac{1}{m}} \right)^m \right]^2,\\
  \label{vGs}
  s(u) = \left \{ 1+\left [ \alpha (u_g-u ) \right]^{\frac{1}{1-m}} \right \}^{-m}.
\end{gather}
Here, $K_s$ is the saturated hydraulic conductivity, $u_g = \frac{P_g}{\rho g}$,  $P_g$ is the air pressure, $\rho$ is the density, $g$ is gravity, and $\alpha$ and $m $ are the van Genuchten parameters. The following parameter values are used in the STOMP simulation: $u_0 = -10$ m, $\alpha = 0.1$, $m = 0.469$, $q=8.25 \times 10^{-5}$ m/s, $u_g = 0$, and $K_s = 8.25 \times 10^{-4}$ m/s.

\begin{figure}[htbp]
  \centering
  \includegraphics[width=0.75\textwidth]{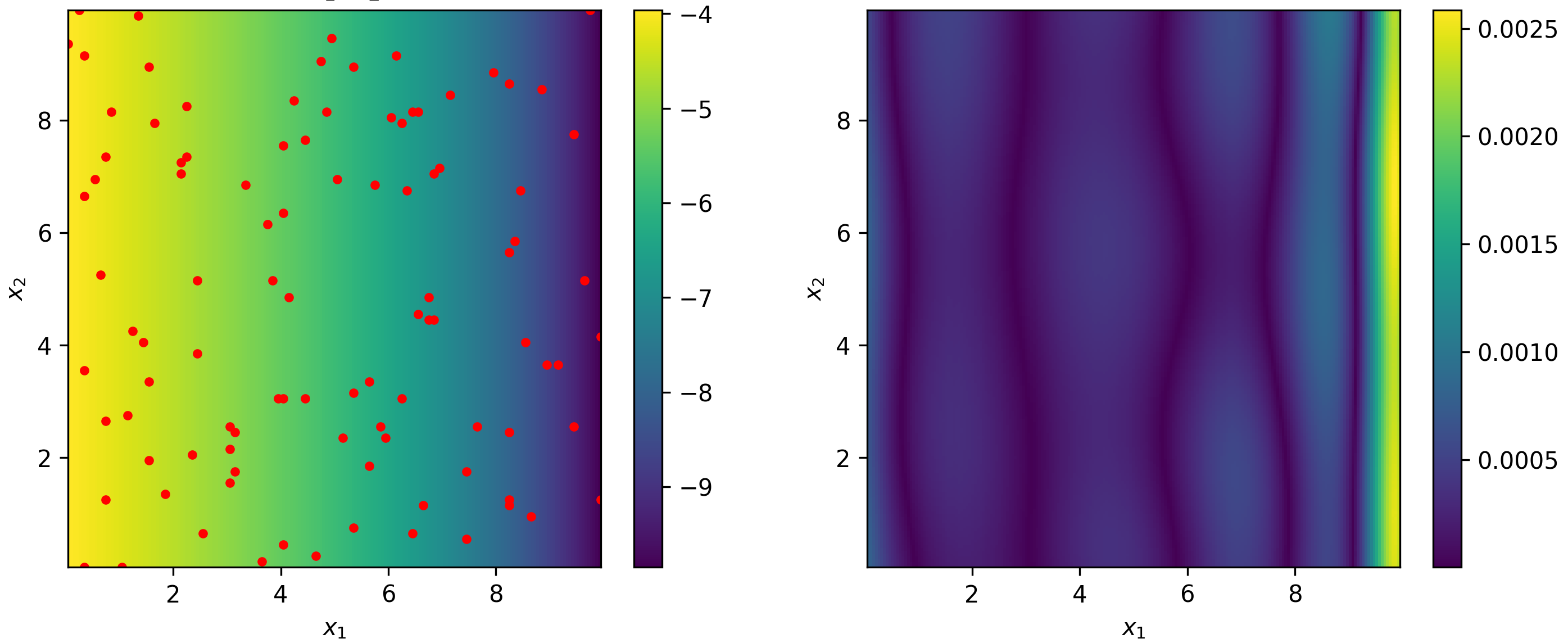}
  \caption{(Left) Referenced $u(\mathbf{x})$ field generated using STOMP with the van Genuchten model for $K(u)$ and the locations of $u$ observations. (Right) Absolute error in the  $u(\mathbf{x})$ field estimated with the physics informed DNN.}
  \label{fig:unsaturated}
\end{figure}

\begin{figure}[htbp]
  \centering
  \includegraphics[width=0.75\textwidth]{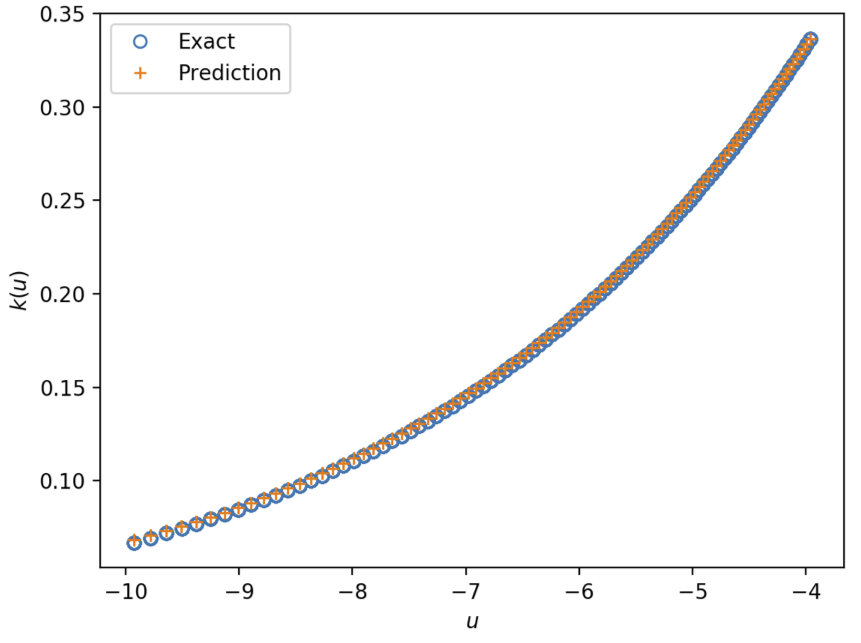}
  \caption{Comparison of the estimated $\hat{K}(u)$ and the reference $K(u)$ given by the van Genuchten model.}
  \label{fig:vanGenuchten}
\end{figure}

Figure \ref{fig:unsaturated} shows the reference $u(\mathbf{x})$ field generated with STOMP and the assumed locations of $u$ measurements. Figure \ref{fig:vanGenuchten} presents the estimated $\hat{K}(u)$ function and the reference $K(u)$ function given by Eqs  (\ref{vGK}) and (\ref{vGs}). It is evident that the physics informed DNN method provides an accurate estimate of unknown $K(u)$ without any direct measurements of $K$ as a function of $u$.

\begin{figure}[htbp]
  \centering
  \includegraphics[width=0.4\textwidth]{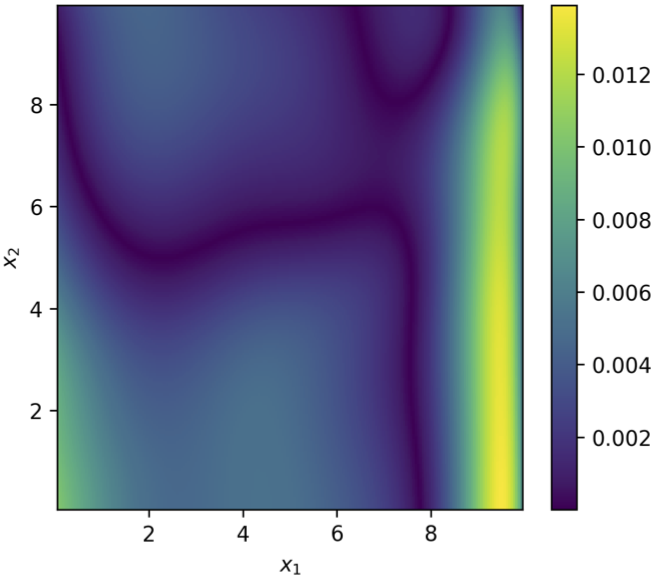}
   \includegraphics[width=0.5\textwidth]{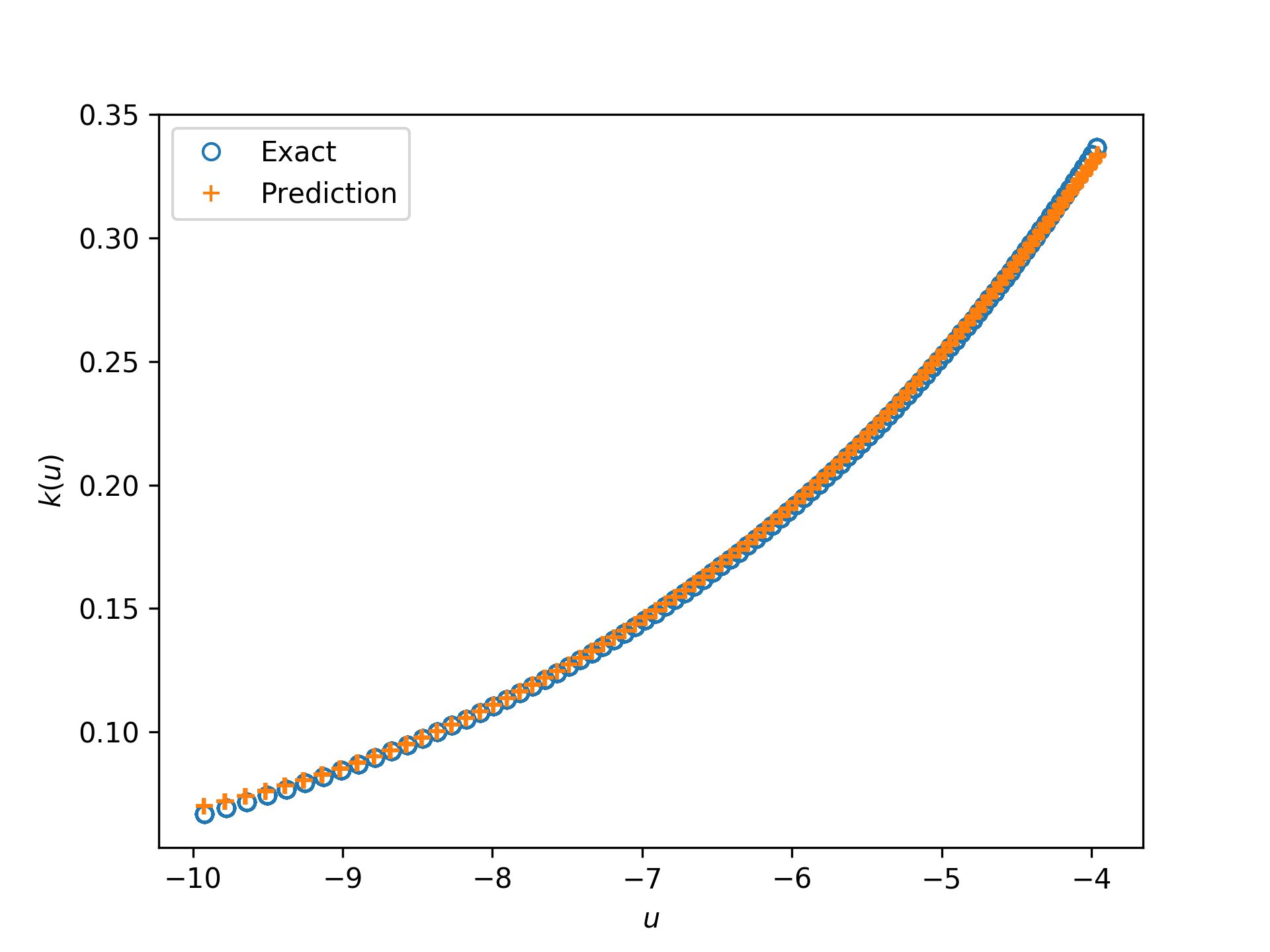}
   \caption{ (Left) Absolute error in $\hat{u}(\mathbf{x})$ in the presence of measurement noise. (Right) Comparison of the reference $K(u)$ given by the van Genuchten model and estimated $\hat{K}(u)$ in the presence of measurement noise.
     }
  \label{fig:unsaturated_noise}
\end{figure}

Finally, we examine the robustness of the physics informed DNN in the presence of observation noise. We use the exact same setup as before except we add 1\% random noise to the values of observed $u$. Figure \ref{fig:unsaturated_noise} shows the difference between predicted and referenced $u(x)$ and $K(u)$.  The added noise increases maximum error in the reconstructed $u$ from 0.002 to 0.012, but the accuracy of the reconstructed $K(u)$ practically does not change. Note that the relative $L_2$ prediction errors for the ``noisy'' case are quite small, including $7.4\times10^{-4}$ for $u$ and $6.4\times10^{-3}$ for $K$. For comparison, the relative $L_2$ prediction errors for the ``noiseless'' case are $5.8\times10^{-5}$ for $u$ and $5.9\times10^{-3}$ for $K$. 

\section{Conclusions}\label{conclusions}

In this work, we have presented a physics informed DNN method for estimating parameters and unknown physics (constitutive relationships) in PDE models. The proposed method uses both PDEs and measurements to train DNNs to approximate unknown parameters and constitutive relationships, as well as states (the PDE solution). Physical knowledge increases the accuracy of DNN training with small data sets and affords the ability to train DNNs when no direct measurements of the functions of interest are available.

We have tested this method for estimating an unknown space-dependent diffusion coefficient in a linear diffusion equation and an unknown constitutive relationship in a non-linear diffusion equation. For the parameter estimation problem, we assume that partial measurements of the coefficient and state are available and have demonstrated that the proposed method is more accurate than the state-of-the-art MAP parameter estimation method. For the non-linear diffusion PDE model with unknown constitutive relationship (state-dependent diffusion coefficient), the proposed method has proven able to accurately estimate the non-linear diffusion coefficient without any measurements of the diffusion coefficient and with measurements of the state only. 
We also have demonstrated that adding physics constraints in training DNNs could increase the accuracy of DNN parameter estimation by as much as 50\%. 

Parameter estimation is an ill-posed problem, and standard parameter estimation methods, including MAP, rely on regularization. In this work, we have trained DNNs for unknown parameters without regularizing  the estimated parameter field or unknown function. In the absence of regularization, we have found that the estimates of the parameter and state mildly depend on the DNN Xavier initialization scheme. For the considered problem, the uncertainty (standard deviation) and mean error decreased with an increasing number of parameter measurements. The coefficient of variation of the relative error (the ratio of the relative error standard deviation to the mean value) was found to be approximately $0.1$. 
In future research, we will investigate the effect of regularization in the physics informed DNN method on  parameter estimation accuracy.

\section{Acknowledgments}
This research was partially supported by the U.S. Department of Energy (DOE) Advanced Scientific Computing (ASCR) and Biological Environmental Research (BER) offices and the  Pacific
Northwest National Laboratory (PNNL) “Deep Learning for Scientific Discovery Agile Investment program. PNNL is operated by Battelle for the DOE under Contract DE-AC05-76RL01830.

%\section*{References}
\bibliographystyle{plain}
%\bibliography{ref.bib}

\end{document}